\documentclass[10pt]{article}
\title{Beyond Monotone Operators: Zarantonello-type Surjectivity Theorems and the Proximal Point Algorithm}
\author{Jos\'e M. M. Brito\thanks{School of Science, Great Bay University, Dongguan 523000, Guangdong Province, People's Republic of China.
E-mail: marciomb@ufpi.edu.br. ORCID-ID: 0000-0003-4362-0536}
\and
Felipe Lara\thanks{Instituto de Alta Investigaci\'on, Universidad de Tarapac\'a, Arica, Chile. E-mail: flarao@academicos.uta.cl.
ORCID-ID: 0000-0002-9965-0921}
\and
Jinyun Yuan\thanks{School of Science, Great Bay University, Dongguan 523000, Guangdong Province, People's Republic of China;
and School of Computer Science and Technology, Dongguan University of Technology,
Dongguan 523808, Guangdong Province, People's Republic of China.
Corresponding author. E-mail: yuanjy570522@126.com. ORCID-ID: 0000-0003-0358-2862}
}

\usepackage[T1]{fontenc}
\usepackage[utf8]{inputenc}
\usepackage[english]{babel}
\usepackage{amsmath,amssymb,amsthm,mathtools}
\usepackage{enumitem}
\usepackage{algorithm}
\usepackage{graphicx}
\usepackage{xcolor}
\usepackage{comment}
\usepackage[hidelinks]{hyperref}
\definecolor{pendingreference}{RGB}{0,70,190}

\newtheorem{theorem}{Theorem}[section]
\newtheorem{proposition}[theorem]{Proposition}
\newtheorem{lemma}[theorem]{Lemma}
\newtheorem{corollary}[theorem]{Corollary}
\theoremstyle{definition}
\newtheorem{definition}[theorem]{Definition}
\newtheorem{example}[theorem]{Example}
\theoremstyle{remark}
\newtheorem{remark}[theorem]{Remark}

\newcommand{\ip}[2]{\left\langle #1,#2\right\rangle}
\newcommand{\norm}[1]{\left\lVert #1\right\rVert}
\newcommand{\Id}{\operatorname{Id}}
\newcommand{\dom}{\operatorname{dom}}

\newcommand{\zer}{\operatorname{zer}}
\newcommand{\gph}{\operatorname{gph}}
\newcommand{\ran}{\operatorname{ran}}
\newcommand{\argmin}{\operatorname{argmin}}
\newcommand{\argmax}{\operatorname{argmax}}
\newcommand{\R}{\mathbb{R}}
\newcommand{\Sph}{\mathbb{S}}

\begin{document}
\maketitle

\begin{abstract}
\noindent 

Outside the maximally monotone class, proximal point iterates need not exist, and their existence is usually assumed. We study the existence and convergence of inexact proximal iterates for set-valued star-monotone operators \(T\) on real Hilbert spaces, with \(\langle u,x-x^\ast\rangle\geq0\) for some zero \(x^\ast\) and every \((x,u)\) in the graph of \(T\). This class strictly contains monotone and pseudomonotone operators with zeros. We prove that \(\mathrm{Id}+cT\) is surjective for every \(c>0\) under four pairwise incomparable topological frameworks, several of whose hypotheses are necessary. Surjectivity can fail under the weak Minty condition with positive violation. For star \(\gamma\)-hypomonotone operators, \(0<c<1/\gamma\) is the sharp range for guaranteed surjectivity. When \(T\) is star-monotone at every zero, we give conditions for weak convergence of inexact sequences with summable errors, allowing arbitrary selections from the possibly set-valued resolvent. A counterexample exhibits nonconvergence under star monotonicity at a single zero. Under strong star monotonicity, we obtain strong convergence, including for certain nonsummable errors, and derive explicit R-linear rates. For a class of nonsmooth nonconvex–nonconcave minimax problems with star quasiconvex partial functions, every proximal step exists for every \(c>0\). In the strongly star quasiconvex case, we obtain R-linear convergence results with subproblems solved to Clarke stationarity.



\medskip

\end{abstract}

\section{Introduction} \label{sec:01}

Let $H$ be a real Hilbert space and let $T:H\rightrightarrows H$ satisfy $\zer(T)\neq\emptyset$. Given $x_0\in H$ and positive parameters $\{c_k\}_{k\geq0}$, we consider the inexact proximal point iteration
\begin{equation}\label{ppa:intro}
 x_k+ e_k = x_{k+1}+ c_k\xi_{k+1}, \qquad \xi_{k+1} \in T(x_{k+1}), \qquad \norm{e_k} \leq \varepsilon_k,
\end{equation}
that is, $x_{k+1}\in J_{c_k}(x_k+e_k)$, where $J_c:=(\Id+cT)^{-1}$. Zero errors recover the exact proximal point algorithm \cite{M1,M2,rock-1976}. For maximally monotone operators, Minty's theorem \cite{Minty1962} gives $\ran(\Id+cT)=H$, and $J_c$ is single-valued and firmly nonexpansive \cite{C1,rock-1976}. Outside this class, two issues arise: the iterates in \eqref{ppa:intro} may fail to exist, and the classical convergence analysis, based on nonexpansiveness of $J_c$, is no longer available.

Our motivation comes from nonsmooth nonconvex--nonconcave minimax problems
\[
 \min_{x\in H_1}\max_{y\in H_2}\ \mathcal L(x,y),\qquad \mathcal L(x,y)=\phi_0(x)-\psi_0(y)-\ip{Ax}{y},
\]
with $\phi_0,\psi_0$ locally Lipschitz and possibly nonconvex, and $A$ bounded and linear. The natural operator is the Clarke saddle operator $T_{\mathcal L}$, the nonsmooth counterpart of Rockafellar's saddle operator \cite{Rockafellar1970}. Each proximal step is itself a nonconvex--nonconcave subproblem, solved in practice only to Clarke stationarity, so the resolvent of $T_{\mathcal L}$ is set-valued. Typical partial functions are quasiconvex functions and their separable sums, which need not be quasiconvex: Bergson--Samuelson social welfare functions in economics \cite{Lara-10}, the Cauchy loss of robust regression \cite{Khanh-Lara-2}, and capped losses of robust representation learning (Example~\ref{ex:robust-alignment}). For such problems, $T_{\mathcal L}$ is not monotone, need not be Lipschitz, and need not be hypomonotone for any modulus (Example~\ref{ex:nonconvex-game}), so neither the existence of the proximal steps nor the convergence of the iterates follows from the available theory. However, whenever the partial functions are star quasiconvex at a saddle point $a$ \cite{Khanh-Lara-1}, the operator $T=T_{\mathcal L}$ satisfies, with $x^\ast=a$, the single radial inequality (Section~\ref{sec:06})
\begin{equation}\label{star:intro}
 \ip{u}{x-x^\ast}\geq0, ~ \forall ~ x\in H, ~ \forall ~ u\in T(x).
\end{equation}
This paper shows that \eqref{star:intro}, together with mild topological assumptions, suffices for both the existence and the convergence of the proximal point algorithm.

Condition \eqref{star:intro}, with $x^\ast\in\zer(T)$, is called star monotonicity at $x^\ast$ in \cite{Khanh-Lara-2} (see also \cite{NTV}); as an algorithmic assumption, it goes back (at least) to \cite{SS}. It compares each graph point only with the zero $x^\ast$. It holds for monotone and pseudomonotone operators with a zero and for Clarke subdifferentials of locally Lipschitz quasar-convex \cite{ABKLY} and star quasiconvex functions \cite{Khanh-Lara-1,Khanh-Lara-2}; for $T=F+N_C$, where $N_C$ is the normal cone to a nonempty closed convex set $C\subset H$, it states that $x^\ast$ is a weak (or Minty) solution of the variational inequality associated with $F$ and $C$ \cite{LJ2025,SZZJM2020}. It is the zero-parameter case of the weak Minty condition \cite{Bohm,DDJ,GTHG,PMC}, it is called star-monotone around $x^\ast$ in \cite{GLG2022}, and, when imposed at every zero, restricted monotonicity in \cite{BBG2022}.

Proximal methods beyond monotonicity have been developed under (co)hy\-po\-mo\-no\-to\-ni\-city-type assumptions \cite{CP,IPS2003,Pennanen2002}, where the analysis is local and well-definedness requires separate justification (see also \cite{ReemReich}), under submonotonicity and metric subregularity \cite{LT-2025}, and under oblique weak Minty conditions in $\R^n$ \cite{EPLP2025}. In \cite{EPLP2025,LT-2025}, the existence of the resolvent steps is a standing hypothesis; in nonconvex--nonconcave min-max optimization, the proximal step is guaranteed only for small parameters or under maximality in a generalized monotonicity class \cite{GTHG,PMC}; and \cite{BBG2022} works in finite dimensions without errors. Explicit methods under \eqref{star:intro}, such as extragradient-type and optimistic dual extrapolation schemes \cite{GLG2022,SZZJM2020}, assume single-valuedness and Lipschitz continuity. Moreover, strongly star-monotone operators need not be hypomonotone, cohypomonotone, nor submonotone near their zero (Example~\ref{ex:not-hadamard}).

Surjectivity theory has advanced by successive weakenings of monotonicity: Zarantonello's contraction and Schauder arguments \cite{Zaran}, his numerical-range criterion \cite{Zaran-2}, Minty's theorem, the Browder--Minty theory \cite{Browder1965, Browder1965-2}, operators of monotone type and Condition~$(M)$ \cite{Brezis1968, BrowderHess1972}, generalized-monotone variational inequalities \cite{BHS2004, HKS, KM}, and maximally comonotone operators \cite{BMW}. Each step relies on pairwise information that \eqref{star:intro} does not provide, and maximality is uninformative for star-monotone operators (see Remark \ref{rem:maximality}). What survives is the topological route of Leray--Schauder type \cite{Park,Zeidler}, which already underlies Zarantonello's global result \cite[Theorem~3]{Zaran}.

Our tools are classical: fixed point theorems of Leray--Schauder type for existence, and quasi-Fejér monotonicity for convergence. The point is that they still work under \eqref{star:intro} alone, where both the pairwise structure behind Minty's theorem and the nonexpansiveness of $J_c$ are lost. To our knowledge, three features are new: the existence of every proximal step for every $c>0$ without maximality; convergence for every selection of a set-valued resolvent, based on a radial substitute for firm nonexpansiveness (Lemma~\ref{prop:pointwise}); and the sharpness of the results, since the hypotheses are tested against examples showing that they cannot be weakened in the directions one would first try.

\medskip

{\bf Our contributions:}
\begin{enumerate}[label=\textnormal{(\roman*)},leftmargin=2.2em]
 \item \emph{Existence.} If $T$ is star-monotone, then $\ran(\Id+cT)=H$ for every $c>0$ under any of four pairwise incomparable topological frameworks, one of which contains the full-domain maximally monotone operators with a zero (Theorem~\ref{thm:surj-unified}). Several of these hypotheses cannot be dropped, and surjectivity can fail under the weak Minty condition with a positive violation parameter. For star $\gamma$-hypomonotone operators, the same frameworks guarantee surjectivity exactly for $0<c<1/\gamma$.

 \item \emph{Convergence for every selection.} If $T$ is star-monotone at every zero, every inexact proximal sequence with summable errors converges weakly to a zero, whichever point of the possibly set-valued resolvent is selected. Unlike \cite{BBG2022,EPLP2025}, the steps are proved to exist, errors are allowed, and $H$ is a Hilbert space. A counterexample on $\R^2$ shows that star monotonicity at a single zero is not sufficient (Example~\ref{ex:failure}). Under strong star monotonicity, convergence is strong, also for nonsummable errors, with explicit rates (Theorem~\ref{thm:strong-convergence}).

 \item \emph{Minimax problems.} For the minimax problems above, whenever the partial functions are star quasiconvex at a saddle point (a class that also contains strongly quasar-convex functions), every proximal step exists for every $c>0$ (Theorem~\ref{thm:saddle-ppa}). In the strongly star quasiconvex case, the iterates converge R-linearly to the saddle point even when the subproblems are solved only to Clarke stationarity, and a gap estimate holds in the quasar-convex case (Corollaries \ref{cor:saddle-sqc} and~\ref{cor:saddle-quasar}).
\end{enumerate}

The paper is organized as follows. Section \ref{sec:02} collects preliminaries. Section \ref{sec:03} proves the Zarantonello-type theorems and their sharpness. Section \ref{sec:05} analyzes the inexact proximal point algorithm, Section~\ref{sec:06} applies the results to nonsmooth minimax problems, and Section~\ref{sec:07} concludes.

\section{Preliminaries and Basic Definitions} \label{sec:02}

Throughout the paper, $H$ is a real Hilbert space with inner product $\ip{\cdot}{\cdot}$ and induced norm $\norm{\cdot}$; strong and weak convergence are denoted by $\to$ and $\rightharpoonup$. We write $B(x,r)$ and $\overline B(x,r)$ for the open and closed balls, and $\overline A$ and ${\rm co}(A)$ for the norm closure and the convex hull of $A\subset H$. Closed balls are weakly compact, and weak compactness is equivalent to weak sequential compactness \cite[Facts 2.34 and 2.37]{bauschke-combettes}. For $T:H\rightrightarrows H$, we use $\dom(T)$, $\ran(T)$, $\gph(T)$ and $\zer(T):=\{x\in H: \, 0\in T(x)\}$ as usual.

\begin{definition} \label{def:semicont}
Let $T:H\rightrightarrows H$ and $\bar x\in H$. Then $T$ is said to be:
\begin{enumerate}[label=\textnormal{(\roman*)},leftmargin=2.2em]
\item \label{def:semicont-norm-lsc} \emph{lower semicontinuous at $\bar x$} if, for every norm-open $V$ with $T(\bar x)\cap V\neq\emptyset$, there is $\delta>0$ such that $T(x)\cap V\neq\emptyset$ for all $x\in B(\bar x,\delta)$;
\item \label{def:semicont-weak-norm-lsc} \emph{weak-to-norm lower semicontinuous at $\bar x$} if the same holds with $B(\bar x,\delta)$ replaced by a weak neighborhood of $\bar x$;
\item \label{def:semicont-norm-usc} \emph{upper semicontinuous at $\bar x$} if, for every norm-open $V\supset T(\bar x)$, there is $\delta>0$ such that $T(x)\subset V$ for all $x\in B(\bar x,\delta)$;
\item \label{def:semicont-weak-usc} \emph{weakly upper semicontinuous at $\bar x$} if, for every weakly open $V\supset T(\bar x)$, there is a weak neighborhood $U$ of $\bar x$ such that $T(x)\subset V$ for all $x\in U$;
\item \label{def:semicont-norm-weak-usc} \emph{norm-to-weak upper semicontinuous at $\bar x$} if \hyperref[def:semicont-weak-usc]{\textnormal{(iv)}} holds with $U$ replaced by a ball $B(\bar x,\delta)$.
\end{enumerate}
The global notions mean that the property holds at every point of $H$.
\end{definition}

These are the standard notions \cite[Proposition 2.1]{michael1956}, \cite[Section 2.5]{BI-book}. Unless otherwise specified, lower and upper semicontinuity refer to items \hyperref[def:semicont-norm-lsc]{\textnormal{(i)}} and \hyperref[def:semicont-norm-usc]{\textnormal{(iii)}}; we write ``norm'' rather than ``strong'' to avoid confusion with strong star monotonicity. An operator $T$ is \emph{compact on bounded sets} if $T(B):=\bigcup_{x\in B}T(x)$ is relatively norm compact for every bounded $B\subset H$; if, in addition, $T$ has closed values, then its values are compact. We use the following elementary facts.

\begin{lemma}[{\cite[Propositions 2 and 3, Theorem 1]{AubinCellina}}] \label{lem:usc-closed-compact}
Let $T:H\rightrightarrows H$. Then the following assertions hold:
\begin{enumerate}[label=\textnormal{(\roman*)},leftmargin=2.2em]
 \item \label{lem:usc-closed} If $T$ has nonempty closed values and is compact on bounded sets, then $T$ is norm upper semicontinuous if and only if $\gph(T)$ is norm closed.
 \item \label{lem:wuc-norm} If $T$ has nonempty weakly compact values, is norm-to-weak upper semicontinuous, and $K\subset H$ is norm compact, then $T(K)$ is weakly compact.
 \item \label{lem:wuc-weak} If $T$ has nonempty weakly compact values, is weakly upper semicontinuous, and $K\subset H$ is weakly compact, then $T(K)$ is weakly compact and $\gph(T)$ is closed in the product of the weak topologies.
\end{enumerate}
\end{lemma}

\begin{definition}[Condition $(M)$] \label{def:condition-M}
 An operator $T:H\rightrightarrows H$ satisfies Condition $(M)$ if, for all nets $\{x_\alpha\}_\alpha$ and $\{u_\alpha\}_\alpha$ in $H$,
\[
 x_\alpha\rightharpoonup x, ~~ u_\alpha\rightharpoonup u, ~~ u_\alpha\in T(x_\alpha), ~~ \limsup_\alpha\ip{u_\alpha}{x_\alpha}\leq\ip{u}{x} ~ \Longrightarrow ~ u\in T(x).
\]
\end{definition}

Condition $(M)$ is the weak-limit closure requirement in Br\'ezis' operators of type $M$ \cite[Definition E]{Brezis1968} (see also \cite{BrowderHess1972,Zeidler}), here formulated for nets.

\begin{remark} \label{rem:M-closed-graph}
Condition $(M)$ implies that $\gph(T)$ is norm closed. If $H$ is finite-dimensional, the converse holds; in infinite dimensions it fails (see $S_3$ in Example~\ref{ex:witnesses}).
\end{remark}

Given $c>0$, the resolvent of $T$ is
\begin{equation} \label{eq:resolvent}
 J_c:=(\Id+cT)^{-1}, \qquad J_c(p):=\{x\in H: \, p\in x+cT(x)\}.
\end{equation}
For $T=\partial f$ with $f$ proper, lower semicontinuous and convex, $J_c$ is the classical proximal mapping. We use three classical tools, where a single-valued mapping is called compact if it maps bounded sets into relatively norm compact sets.

\begin{theorem}[Michael's selection theorem {\cite[Theorem $3.2^{\prime\prime}$]{michael1956}}] \label{thm:michael}
If $X$ is paracompact Hausdorff and $Y$ is Banach, every lower semicontinuous $\Phi:X\rightrightarrows Y$ with nonempty closed convex values admits a continuous selection.
\end{theorem}

\begin{theorem}[Schaefer's fixed-point theorem {\cite{schaefer1955}}] \label{thm:schaefer}
 If $f$ is a continuous compact self-map of a Banach space $X$ and $\{x\in X: \, x=\tau f(x),\ \tau \in [0, 1]\}$ is bounded, then $f$ has a fixed point.
\end{theorem}

\begin{theorem}[Set-valued Leray--Schauder principle {\cite[Theorem 3.1]{avramescu2004}}] \label{thm:ls-set}
If $\Phi: X \rightrightarrows X$ is upper semicontinuous and compact on bounded sets, with nonempty closed convex values, on a Banach space $X$, and $\{x\in X: \, x\in\tau \Phi(x),\ \tau\in \, ]0, 1[\}$ is bounded, then $\Phi$ has a fixed point.
\end{theorem}

An operator $T:K\rightrightarrows H$, with $K\subset H$, is \emph{monotone} on $K$ if $\ip{u-v}{x-y} \geq 0$ and \emph{pseudomonotone} on $K$, in the sense of Karamardian \cite{Karamardian1976}, if $\ip{u}{y-x}\geq0$ implies $\ip{v}{y-x} \geq 0$, for all $x,y\in K$, $u\in T(x)$ and $v\in T(y)$ (see \cite{BI-book,HKS}). The following notion has been studied recently in \cite{Khanh-Lara-2} (see also \cite{NTV}).

\begin{definition}[Star monotonicity] \label{def:star}
Let $T:H\rightrightarrows H$ and $x^\ast\in H$.
Then $T$ is said to be \emph{(strongly) star-monotone at $x^\ast$ with modulus $\gamma\geq 0$} if $0\in T(x^\ast)$ and
\begin{equation}\label{eq:strong-star}
 \ip{u}{x-x^\ast} \geq \gamma \norm{x-x^\ast}^2, ~ \forall ~ x \in H, ~ \forall ~ u \in T(x).
\end{equation}
When $\gamma=0$ in \eqref{eq:strong-star}, we simply say that $T$ is star-monotone at $x^{*}$.
\end{definition}

Every monotone or pseudomonotone operator with a zero $x^\ast$ is star-monotone at $x^\ast$ \cite[page 765]{SS}, and for single-valued monotone hemicontinuous operators, Minty's lemma \cite[Chapter III, Lemma 1.5]{KinderlehrerStampacchia} shows that \eqref{eq:strong-star} with $\gamma=0$ says no more than $T(x^\ast)=0$. Outside the monotone class it is strictly stronger: $T(t):=-t$ on $\R$ has the zero $0$, but $tT(t)<0$ for $t\neq0$ and $\ran(\Id+T)=\{0\}$. Star monotonicity and all the semicontinuity, closedness, convexity and compactness properties above are invariant under the translation $T\mapsto T(\cdot+x^\ast)$ \cite[Proposition 9]{Khanh-Lara-2}, and strong star monotonicity with $\gamma>0$ gives $\zer(T)=\{x^\ast\}$.

We also use the following quadratic relaxations of monotonicity and of star monotonicity.

\begin{definition} \label{def:hypo}
Let $T:H\rightrightarrows H$, $x^\ast\in H$, $\gamma>0$ and $\rho>0$.
Then $T$ is said to be:
\begin{enumerate}[label=\textnormal{(\roman*)},leftmargin=2.2em]
\item $\gamma$-hypomonotone if
$\ip{u-v}{x-y}\geq-\gamma\norm{x-y}^2$ for all $(x,u),(y,v)\in\gph(T)$;
\item $\rho$-cohypomonotone if
$\ip{u-v}{x-y}\geq-\rho\norm{u-v}^2$ for all $(x,u),(y,v)\in\gph(T)$;
 \item star $\gamma$-hypomonotone at $x^\ast$ if
$\ip{u}{x-x^\ast}\geq-\gamma\norm{x-x^\ast}^2$ for all $(x,u)\in\gph(T)$;
 \item star $\rho$-cohypomonotone at $x^\ast$ if
$\ip{u}{x-x^\ast}\geq-\rho\norm{u}^2$ for all $(x,u)\in\gph(T)$.
\end{enumerate}
In \textnormal{(i)} and \textnormal{(ii)}, $T$ is maximal if its graph admits no proper extension within the same class.
\end{definition}

Thus $T$ is $\gamma$-hypomonotone if and only if $T+\gamma\Id$ is monotone \cite{IPS2003, Pennanen2002}, and $\rho$-cohypomonotone if and only if $T^{-1}+\rho\Id$ is monotone \cite{CP}; in the terminology of \cite{BMW}, these are $(-\gamma)$-monotonicity and $(-\rho)$-comonotonicity, respectively. If $0\in T(x^\ast)$, star $\rho$-cohypomonotonicity at $x^\ast$ means that $x^\ast$ is a weak Minty solution with parameter $2\rho$ in the sense of \cite[Assumption 1]{Bohm}; it is also called star-negative comonotonicity in \cite{GTHG}.

\begin{remark}[Maximality is uninformative] \label{rem:maximality}
Every operator star-monotone at $x^\ast$ is contained in $\mathcal H_{x^\ast}(x) := \{u \in H: \, \ip{u}{x-x^\ast}\geq0\}$, which is itself star-monotone at $x^\ast$; hence $\mathcal H_{x^\ast}$ is the unique maximal element of the class. This maximality gives surjectivity only trivially, since $x^\ast + c \mathcal H_{x^\ast}(x^\ast) = H$, and it is not necessary: $T(x)=x-x^\ast$ is not maximal, yet $\Id+cT$ is a bijection. We therefore seek topological sufficient conditions for surjectivity.
\end{remark}

\section{Zarantonello-type Theorems} \label{sec:03}

\subsection{Surjectivity Beyond Monotone Operators}

The following theorem is our main result. It establishes $\ran(\Id+cT)=H$ for every $c>0$ under star monotonicity and any one of four alternative sets of topological assumptions.

\begin{theorem}[Zarantonello-type theorem] \label{thm:surj-unified}
 Let $T:H\rightrightarrows H$ be star-monotone at some $x^\ast\in H$. Assume that at least one of the following conditions holds:
\begin{enumerate}[label=\textnormal{(Z\arabic*)},leftmargin=2.8em]
 \item \label{Z1} $T$ is norm lower semicontinuous, compact on bounded sets, and has non\-empty closed convex values;

 \item \label{Z2} $T$ is weak-to-norm lower semicontinuous and has nonempty closed convex values;
 
 \item \label{Z3} $T$ is norm upper semicontinuous, compact on bounded sets, and has nonempty closed convex values; 
 
 \item \label{Z4} $T$ is norm-to-weak upper semicontinuous, has nonempty weakly compact convex values, and satisfies Condition $(M)$.
\end{enumerate}
Then
\[
 \ran(\Id+cT)=H, ~ \forall ~ c>0.
\]
Equivalently, $J_c(p) \neq \emptyset$ for every $p \in H$ and every $c>0$.
\end{theorem}

\begin{proof}
Fix $p\in H$ and $c>0$, set $q:=p-x^\ast$, and define
\begin{equation*}
 \Phi(w) := q-cT(x^\ast+w), \qquad w\in H.
\end{equation*}
Then, $p\in y+cT(y)$ with $y:=x^\ast+w$ if and only if $w\in\Phi(w)$. Cases \textnormal{(Z1)}--\textnormal{(Z3)} rest on the a priori bound
\begin{equation}\label{eq:apriori}
 w=\tau(q-cu),\quad \tau\in[0,1],\quad u\in T(x^\ast+w) ~ \Longrightarrow ~ \norm{w}\leq\tau\norm{q},
\end{equation}
which holds because star monotonicity gives $\ip{u}{w}\geq0$, so $\norm{w}^2=\tau\ip{q}{w}-\tau c\ip{u}{w}\leq\tau\norm{q}\norm{w}$.

\medskip
\noindent\textnormal{(Z1)}: $\Phi$ is lower semicontinuous with nonempty closed convex values. By Theorem~\ref{thm:michael}, it admits a continuous selection $s:H\to H$, which is compact because 
$$s(B) \subset q-cT(x^\ast+B)$$ 
for bounded $B \subset H$. By \eqref{eq:apriori}, we have $\{w\in H: w=\tau s(w),\ \tau\in[0,1]\} \subset \overline B(0,\norm{q})$, and Theorem \ref{thm:schaefer} yields $\bar w=s(\bar w) \in \Phi(\bar w)$.

\medskip
\noindent\textnormal{(Z2)}: Choose $R>\norm{q}$ and let $K:=\overline B(0,R)$, endowed with the weak topology; $K$ is compact \cite[Fact 2.34]{bauschke-combettes}, hence paracompact. The restriction $\Phi:K\rightrightarrows H$ is lower semicontinuous from the weak topology of $K$ into the norm topology of $H$, with nonempty closed convex values, so Theorem \ref{thm:michael} gives a weak-to-norm continuous selection $s:K\to H$; in particular, $s(K)$ is norm compact. With the metric projection $P_Kz=\dfrac{Rz}{\max\{R,\norm{z}\}}$, the map $f := s\circ P_K$ is norm continuous and compact. 

Let $z=\tau f(z)$ with $\tau\in[0,1]$ and $w:=P_Kz$. If $\norm{z}\leq R$, \eqref{eq:apriori} gives $\norm{z}\leq\norm{q}$. If $\norm{z}>R$, then $w=\lambda s(w)$ with $\lambda := \tau R/\norm{z}\in[0,1[$, and \eqref{eq:apriori} gives $R=\norm{w} \leq \norm{q}$, a contradiction. Theorem \ref{thm:schaefer} yields $\bar z=f(\bar z)$; then $\norm{\bar z}\leq\norm{q}<R$, so $\bar z=s(\bar z)\in\Phi(\bar z)$.

\medskip
\noindent\textnormal{(Z3)}: $\Phi$ is upper semicontinuous and compact on bounded sets, with nonempty closed convex values. By \eqref{eq:apriori}, we have $\{w\in H: w\in \tau \Phi(w),\ \tau\in\,]0,1[\} \subset \overline B(0,\norm{q})$, and Theorem \ref{thm:ls-set} yields $\bar w \in \Phi(\bar w)$.

\medskip

\noindent\textnormal{(Z4)}:
Define $\widehat T(w) := T(x^\ast+w)$. Translation preserves the value assumptions, norm-to-weak upper semicontinuity and condition $(M)$, since $\ip{u_\alpha}{x^\ast +w_\alpha} = \ip{u_\alpha}{w_\alpha} + \ip{u_\alpha}{x^\ast}$ and the last term converges to $\ip{u}{x^\ast}$ whenever $u_\alpha \rightharpoonup u$. Moreover, $0 \in \widehat T(0)$ and $\ip{u}{w} = \ip{u}{(x^\ast+w) - x^\ast} \geq 0$ for every $u\in \widehat T(w)$. Thus $\widehat T$ is star-monotone at $0$.

Let $\mathcal D$ be the directed family of finite-dimensional subspaces
$E \subset H$ containing $q$, ordered by inclusion. For $E\in\mathcal D$, let
$P_E$ be the orthogonal projection onto $E$ and set
\[
 T_E(w):=P_E\widehat T(w), \qquad w\in E.
\]
The operator $T_E$ is norm upper semicontinuous and has nonempty compact convex values. It is compact on bounded subsets of $E$, since the closure of such a set is compact, its image under $\widehat T$ is weakly compact by Lemma \ref{lem:usc-closed-compact}\hyperref[lem:wuc-norm]{\textnormal{(ii)}},
and $P_E$ maps weakly compact sets into norm-compact subsets of $E$.
Moreover, $0\in T_E(0)$ and $\ip{P_Eu}{w}=\ip{u}{w}\geq0$ for every $w\in E$
and $u\in\widehat T(w)$, so $T_E$ is star-monotone at $0$ in $E$. Thus,
case \textnormal{(Z3)} above, applied in the Hilbert space $E$ to $T_E$,
gives $w_E \in E$ and $u_E \in \widehat T(w_E)$ such that
\begin{equation} \label{eq:galerkin-type-M}
 q=w_E+cP_Eu_E.
\end{equation}

Since $\ip{P_Eu_E}{w_E}=\ip{u_E}{w_E}\geq0$, the computation in \eqref{eq:apriori} gives $\norm{w_E}\leq\norm{q}$.
Passing to a subnet, we may assume that $w_E \rightharpoonup w$. For each
fixed $v\in H$, set $E_v:=\operatorname{span}\{q,v\}\in\mathcal D$; then
$v\in E$ for all sufficiently large $E$ (along the original net, for every
$E\supseteq E_v$), and \eqref{eq:galerkin-type-M} then implies
\[
 \ip{u_E}{v} = \ip{P_Eu_E}{v} = \frac1c\ip{q-w_E}{v}
 \longrightarrow \frac1c\ip{q-w}{v}.
\]
Hence $u_E\rightharpoonup u:=(q-w)/c$. Furthermore, since $w_E\in E$,
\[
 \ip{u_E}{w_E} = \ip{P_Eu_E}{w_E}
 = \frac1c\bigl(\ip{q}{w_E}-\norm{w_E}^2\bigr),
\]
and the weak lower semicontinuity of the norm gives
\[
 \limsup_E\ip{u_E}{w_E}
 \leq \frac1c\bigl(\ip{q}{w}-\norm{w}^2\bigr) = \ip{u}{w}.
\]
Then, Condition $(M)$ yields $u\in\widehat T(w)$. Thus, with
$y:=x^\ast+w$, we conclude that $p=y+cu \in y+cT(y)$.
\end{proof}

\begin{remark} \label{rem:closed-graph-IIA}
 Under the other assumptions of 
case \textnormal{(Z3)} of Theorem~\ref{thm:surj-unified}, norm upper semicontinuity may equivalently be replaced by norm closedness of $\gph(T)$ (see Lemma~\ref{lem:usc-closed-compact}\hyperref[lem:usc-closed]{\textnormal{(i)}}). Compactness on bounded sets is essential for this equivalence, even when $H=\R$ (see Example~\ref{ex:hypotheses}\textnormal{(a)} below).
\end{remark}

As far as we know, the unified theorem is also novel for the class of pseudomonotone operators. Indeed, in the presence of a zero, pseudomonotonicity provides precisely the radial inequality required by star monotonicity.

\begin{corollary}[Pseudomonotone operators] \label{cor:pseudomonotone-surj}
 Let $T: H \rightrightarrows H$ be pseudomonotone and $x^\ast\in\zer(T)$. If one of \hyperref[Z1]{\textnormal{(Z1)}}--\hyperref[Z4]{\textnormal{(Z4)}} holds, then $\ran(\Id+cT)=H$ for every $c>0$.
\end{corollary}

\begin{proof}
 Fix $x\in H$ and $u\in T(x)$. Since $x^\ast\in\zer(T)$, $0 \in T(x^\ast)$. Applying pseudomonotonicity to the points $x^\ast$ and $x$, and to the elements $0 \in T(x^\ast)$ and $u\in T(x)$, gives $\ip{0}{x-x^\ast} = 0$ implies $\ip{u}{x-x^\ast} \geq 0$. Therefore, $T$ is star-monotone at $x^\ast$, and the conclusion follows from Theorem \ref{thm:surj-unified}.
\end{proof}

\begin{corollary}[Weak upper semicontinuity] \label{cor:surj-weak-usc}
 Let $T:H\rightrightarrows H$ be weakly upper semicontinuous with nonempty weakly compact convex values. If $T$ is star-monotone at some $x^\ast\in H$, then $\ran(\Id+cT)=H$ for every $c>0$.
\end{corollary}

\begin{proof}
 Weak upper semicontinuity implies norm-to-weak upper semicontinuity. Then, Lemma \ref{lem:usc-closed-compact}\hyperref[lem:wuc-weak]{\textnormal{(iii)}} gives weak-to-weak closedness of the graph, which implies Condition $(M)$. Hence, the result follows by applying case \textnormal{(Z4)} of Theorem \ref{thm:surj-unified}.
\end{proof}

\begin{corollary}[Full-domain maximally monotone operators] \label{cor:maximal-monotone}
 Let $T: H \rightrightarrows H$ be maximally monotone with $\dom(T)=H$ and $\zer(T) \neq \emptyset$. Then $T$ satisfies the assumptions of 
case \textnormal{(Z4)} of Theorem~\ref{thm:surj-unified}. In particular, $\ran (\Id+cT) = H$ for every $c>0$.
\end{corollary}

\begin{proof}
 By \cite[Corollaries~4.2.6, 4.2.11, and~4.2.13]{BI-book}, the values of $T$ are nonempty, weakly compact, and convex, and $T$ is norm-to-weak upper semicontinuous. To verify Condition $(M)$, let $x_\alpha\rightharpoonup x$, $u_\alpha \rightharpoonup u$, $u_\alpha \in T(x_\alpha)$, and assume $\limsup_\alpha \ip{u_\alpha}{x_\alpha} \leq \ip{u}{x}$. Then, for every $(y,v)\in\gph(T)$, monotonicity gives
\[
 0 \leq\limsup_\alpha\ip{x_\alpha-y}{u_\alpha-v} \leq \ip{x-y}{u-v}.
\]
Hence, $(x, u)$ is monotonically related to $\gph(T)$, and maximality gives $u \in T(x)$. Mo\-no\-to\-ni\-ci\-ty relative to any $x^\ast \in\zer(T)$ implies star monotonicity at $x^\ast$.
\end{proof}

For maximally monotone operators, the zero set is closed and convex \cite[Proposition 23.39]{bauschke-combettes}. The surjectivity in Theorem~\ref{thm:surj-unified} recovers this property beyond monotonicity, although the resolvent may be set-valued and not Lipschitz continuous (see Example~\ref{rem:resolvent-failure}).

\begin{corollary}[Closedness and convexity of the zero set] \label{cor:zero-set}
Let $T: H \rightrightarrows H$ be star-monotone at every point of its nonempty zero set, and suppose that $\ran(\Id+cT)=H$ for some $c>0$, for instance under one of \hyperref[Z1]{\textnormal{(Z1)}}--\hyperref[Z4]{\textnormal{(Z4)}}. Then $\zer(T)$ is closed and convex.
\end{corollary}

\begin{proof}
Set $C := \overline{\mathrm{co}}(\zer(T))$ and let $p \in C$. By surjectivity, there exists $x\in J_c(p)$. Star monotonicity at each $z \in \zer(T)$ gives $\ip{p-x}{x-z} \geq 0$, and by linearity and continuity in $z$ this inequality holds for every $z \in C$. Taking $z=p$ yields $-\norm{p-x}^2 \geq 0$, so $x=p$ and $0 \in T(p)$. Hence $C=\zer(T)$.
\end{proof}

Star monotonicity at every zero is essential: the continuous operator of Example \ref{ex:depends-zero} is star-monotone at $0$ and satisfies \hyperref[Z3]{\textnormal{(Z3)}}, but $\zer(T)=\{0\}\cup\Sph^1$ is not convex. Solvability is also essential: $T(x,y):=(0,y)$ for $y\neq0$ and $T(x,0):=(0,x(x-1))$ is star-monotone at both of its zeros $(0,0)$ and $(1, 0)$, but $J_c(1/2,0)=\emptyset$ for every $c>0$.

\begin{remark}[Finite-dimensional reduction] \label{rem:finite-reduction}
If $H=\R^n$, the weak and norm topologies coincide, and, by Remark~\ref{rem:M-closed-graph}, Theorem~\ref{thm:surj-unified} reduces to two alternatives: $T$ is lower semicontinuous with nonempty closed convex values, or $T$ is upper semicontinuous with nonempty compact convex values.
\end{remark}

\subsection{Comparison and Sharpness} \label{sec:comparison}

The closest classical result, \cite[Theorem 19]{Brezis1968}, requires a single-valued operator of type~$M$, bounded on bounded sets, with $\ip{A(x)}{x}\geq0$ on a sphere; it is not comparable with \hyperref[Z4]{\textnormal{(Z4)}}, which needs no boundedness but imposes the sign condition on all of $H$.


\begin{example}[Separating operators] \label{ex:witnesses}
Let $H$ be infinite-dimensional, $e_1\in H$ be a unit vector, $\{e_k\}_{k\geq2}\subset e_1^\perp$ be orthonormal, and set $\mathcal H := \R \times H$. Define $Q(x):=\norm{x}\ip{x}{e_1}e_1$ on $H$, and on $\R$ define $L(0):=\{0\}$, $L(t) := \operatorname{co}\{t, \operatorname{sgn}(t)\}$ for $t\neq0$, and $U(t) := \{\operatorname{sgn}(t)\}$ for $t\neq0$, $U(0):=[-1,1]$. Consider
\begin{enumerate}[label=\textnormal{(\alph*)},leftmargin=2.2em]
 \item $S_1(t,x):=L(t)\times\{Q(x)\}$ on $\mathcal H$;
 \item $S_2(0):=\{0\}$ and $S_2(x):=\{u\in H:\ip{u}{x}\geq\ip{x}{e_1}^2\}$ for $x\neq0$;
 \item $S_3(t,x):=U(t)\times\{Q(x)\}$ on $\mathcal H$;
 \item $S_4(x):=x+\norm{x}e_1$ on $H$.
\end{enumerate}
All four are star-monotone at the origin. Indeed, every $u\in L(t)\cup U(t)$ satisfies $ut\geq0$, and $\ip{Q(x)}{x}=\norm{x}\ip{x}{e_1}^2\geq0$, which covers $S_1$ and $S_3$; for $S_2$ this is built into the definition; and $\ip{S_4(x)}{x}=\norm{x}^2+\norm{x}\ip{e_1}{x}\geq0$.
\end{example}

\begin{proposition}[Sharp separation of the frameworks] \label{prop:table}
The operators of Example~\ref{ex:witnesses} satisfy the frameworks of Theorem~\ref{thm:surj-unified} exactly as follows:
\begin{center}
\begin{tabular}{lcccc}
\hline
& \hyperref[Z1]{\textnormal{(Z1)}} & \hyperref[Z2]{\textnormal{(Z2)}} & \hyperref[Z3]{\textnormal{(Z3)}} & \hyperref[Z4]{\textnormal{(Z4)}}\\
\hline
$S_1$ & yes & no  & no  & no  \\
$S_2$ & no  & yes & no  & no  \\
$S_3$ & no  & no  & yes & no  \\
$S_4$ & no  & no  & no  & yes \\
\hline
\end{tabular}
\end{center}
In particular, the four frameworks are pairwise incomparable, and each has an exclusive representative.
\end{proposition}

\begin{proof}
Two observations are used repeatedly. First, with $x_k:=e_1+e_k$ we have $x_k\rightharpoonup e_1$ and $Q(x_k)=\sqrt2\,e_1\neq e_1=Q(e_1)$, so $Q$ is norm continuous and compact on bounded sets (its range lies in $\R e_1$), but not weak-to-norm continuous. Second, on the factor $\R$ the weak and norm topologies coincide, so a failure of semicontinuity of $L$ or $U$ at $0$ passes to $S_1$ or $S_3$ by fixing the $H$-component at $0$.

\medskip

\noindent $(a)$: $L$ is lower semicontinuous with compact convex values and maps bounded sets into bounded sets; with the properties of $Q$, $S_1$ satisfies \hyperref[Z1]{\textnormal{(Z1)}}. Since $Q$ is not weak-to-norm continuous, \hyperref[Z2]{\textnormal{(Z2)}} fails. Since $L(t)=[t,1]\not\subset\,]-\tfrac12,\tfrac12[\,\supset L(0)$ for $t\in\,]0,\tfrac12[$, $L$ is not upper semicontinuous at $0$, which excludes \hyperref[Z3]{\textnormal{(Z3)}} and \hyperref[Z4]{\textnormal{(Z4)}}.

\medskip

\noindent $(b)$: The values of $S_2$ are closed and convex, and unbounded for $x\neq0$, which excludes \hyperref[Z1]{\textnormal{(Z1)}}, \hyperref[Z3]{\textnormal{(Z3)}} and \hyperref[Z4]{\textnormal{(Z4)}}. For weak-to-norm lower semicontinuity at $\bar x\neq0$, let $\bar u\in S_2(\bar x)$ and $t>0$: the map $x\mapsto\ip{\bar u+t\bar x}{x}-\ip{x}{e_1}^2$ is weakly continuous and positive at $\bar x$, so $\bar u+t\bar x\in S_2(x)$ on a weak neighborhood of $\bar x$ not containing $0$; letting $t\downarrow0$ gives the claim. At $0$, for $x\neq0$ with $|\ip{x}{e_1}|<\varepsilon$, the vector $u_x:=\ip{x}{e_1}^2\norm{x}^{-2}x$ lies in $S_2(x)$ and $\norm{u_x}<\varepsilon$. Hence \hyperref[Z2]{\textnormal{(Z2)}} holds.

\medskip

\noindent $(c)$: $U$ is upper semicontinuous with compact convex values in $[-1,1]$; with the properties of $Q$, $S_3$ satisfies \hyperref[Z3]{\textnormal{(Z3)}}. Since $U(t)=\{-1\}$ misses $]\tfrac12,\tfrac32[$, which meets $U(0)$, for $t<0$, \hyperref[Z1]{\textnormal{(Z1)}} and \hyperref[Z2]{\textnormal{(Z2)}} fail. Condition $(M)$ fails as well: with $z_k:=(0,x_k)\rightharpoonup z:=(0,e_1)$ and $w_k=w:=(0,\sqrt2\,e_1)\in S_3(z_k)$, we have $\ip{w_k}{z_k}=\sqrt2=\ip{w}{z}$, but $w\notin S_3(z)=[-1,1]\times\{e_1\}$. Hence \hyperref[Z4]{\textnormal{(Z4)}} fails.

\medskip

\noindent $(d)$: $S_4$ is norm continuous and single-valued, hence norm-to-weak upper semicontinuous with weakly compact convex values. For Condition $(M)$, let $x_\alpha \rightharpoonup x$, $S_4(x_\alpha) \rightharpoonup u$ and $\limsup_\alpha \ip{S_4(x_\alpha)}{x_\alpha}\leq\ip{u}{x}$. Since $\norm{x_\alpha}e_1\rightharpoonup u-x$, we get $\norm{x_\alpha}\to\rho:=\ip{u-x}{e_1}\geq0$ and $u=x+\rho e_1$. The inequality then reads $\rho^2+\rho\ip{e_1}{x}\leq \norm{x}^2+\rho\ip{e_1}{x}$, so $\rho\leq \norm{x}$, while weak lower semicontinuity of the norm gives $\norm{x}\leq\rho$. Thus $u=S_4(x)$, and \hyperref[Z4]{\textnormal{(Z4)}} holds. Finally, $\norm{S_4(e_k)-S_4(e_m)} =\sqrt2$ for $k\neq m$, so $S_4$ is not compact on bounded sets, which excludes \hyperref[Z1]{\textnormal{(Z1)}} and \hyperref[Z3]{\textnormal{(Z3)}}; and $e_k\rightharpoonup0$ with $\norm{S_4(e_k)-S_4(0)}=\sqrt2$, which excludes \hyperref[Z2]{\textnormal{(Z2)}}.

Therefore, the four frameworks are pairwise incomparable, and each has an exclusive representative.
\end{proof}

In particular, \hyperref[Z2]{\textnormal{(Z2)}} is the only framework that allows unbounded values. Since Theorem~\ref{thm:surj-unified} is a disjunction, a hypothesis of one framework may be replaced by those of another; the examples below show that several hypotheses cannot simply be dropped (on $\R$, Condition $(M)$ amounts to closedness of the graph; see Remark~\ref{rem:M-closed-graph}).

\begin{remark}[Sharpness of star monotonicity] \label{rem:weak-minty}
Star monotonicity cannot be relaxed to star cohypomonotonicity (Definition~\ref{def:hypo}), that is, to the weak Minty condition of
\cite{Bohm,DDJ} with a positive parameter. Given $\rho>0$ and $c>0$, set $a:=\max\{1/\rho,2/c\}$ and $T(t):=-at$ for $t\geq0$, $T(t):=t$ for
$t<0$. Since $a\rho\geq1$, we have $tT(t)\geq-\rho T(t)^2$ for every $t\in\R$, so $T$ is star $\rho$-cohypomonotone at $0$; moreover, $T$ is
continuous and compact on bounded sets, so it satisfies every topological and value requirement of
\hyperref[Z1]{\textnormal{(Z1)}}--\hyperref[Z4]{\textnormal{(Z4)}}.
However, since $ac>1$,
\[
 \ran(\Id+cT)=\{(1-ac)t: \, t\geq0\}\cup\{(1+c)t: \, t<0\}=\,]-\infty,0]\neq\R.
\]
Thus the value $\rho=0$, that is, star monotonicity, cannot be improved within this family, for any $c>0$. In contrast, star hypomonotonicity preserves surjectivity for small parameters (Corollary~\ref{cor:weak-star}).
\end{remark}

\begin{example}[Role of the hypotheses] \label{ex:hypotheses}
Each of the following operators on $\R$ is star-monotone at $0$ and fails one requirement of the frameworks, and $\Id+cT$ is not surjective. Similar elementary examples show that nonemptiness and convexity of the values cannot be dropped.

\medskip

\noindent $(a)$: \emph{Upper semicontinuity in \textnormal{(Z4)}.} Let $T(t):=0$ for $t\leq0$ and $T(t):=1/t$ for $t>0$. Then $T$ has singleton values and a closed graph, hence satisfies Condition $(M)$, but it is not upper semicontinuous at $0$ (nor compact on bounded sets). Since $t+c/t\geq2\sqrt c$ for $t>0$, we get $\ran(\Id+cT)=\,]-\infty,0]\cup[2\sqrt c,+\infty[$ for every $c>0$. Hence upper semicontinuity cannot be replaced by closedness of the graph, even with compact values.

\medskip

\noindent $(b)$: \emph{Semicontinuity in \textnormal{(Z1)--(Z3)}.} Let $T(t):=\{0\}$ for $t \leq0$ and $T(t):=\{1\}$ for $t>0$. Then $T$ has compact convex values and is compact on bounded sets, but it is neither lower nor upper semicontinuous at $0$, and $\ran(\Id+cT) = \, ]-\infty, 0] \cup \, ]c, + \infty[$ for every $c>0$.
\end{example}

\subsection{Star Hypomonotone and Star Cohypomonotone Operators} \label{subsec:hypo}

We now consider the star hypomonotone operators of Definition \ref{def:hypo}. If $T$ is maximally $\gamma$-hypomonotone, then $\ran(\Id+cT)=H$ for every $c<1/\gamma$ \cite{BMW}; for single-valued operators, the same bound follows from Zarantonello's numerical-range theorem \cite[Theorems II and III]{Zaran-2}. The following result recovers it from the radial inequality alone.

\begin{corollary} \label{cor:weak-star}
 Let $T:H\rightrightarrows H$ be star $\gamma$-hypomonotone at some $x^\ast \in H$, with $\gamma>0$, and assume that one of
\hyperref[Z1]{\textnormal{(Z1)}}--\hyperref[Z4]{\textnormal{(Z4)}} holds. Then $\ran(\Id+cT)=H$ for every $c\in\,]0,1/\gamma[$. Hence, for every $x_0\in H$ and every $\{c_k\}_{k \geq 0} \subset \,]0,1/\gamma[$, there exists a sequence $\{x_k\}_{k \geq 0}$ with $x_{k+1} \in J_{c_k}(x_k)$ for all $k \geq 0$.
\end{corollary}

\begin{proof}
Repeat the proof of Theorem~\ref{thm:surj-unified}. The radial
inequality is used only in \eqref{eq:apriori}, which now reads
$\norm{w}^2\leq\tau\ip{q}{w}+\tau c\gamma\norm{w}^2$, and gives
$(1-c\gamma)\norm{w}\leq\norm{q}$. In case \textnormal{(Z2)}, choose
$R>\norm{q}/(1-c\gamma)$ and use that $\lambda\mapsto
\lambda/(1-\lambda c\gamma)$ is increasing on $[0,1]$; in case
\textnormal{(Z4)}, note that $T_E$ is star $\gamma$-hypomonotone at
$0$ in $E$.
\end{proof}

\begin{remark}[Sharpness of the bound $1/\gamma$]
\label{rem:weak-star-sharp}
 The bound $c<1/\gamma$ in Corollary \ref{cor:weak-star} cannot be improved. Let $T:\R\to\R$ be given by $T(t):=-\gamma t$ for $t\geq0$ and $T(t):=t$ for $t<0$, that is, the operator of Remark \ref{rem:weak-minty} with $a=\gamma$. Then $T(0)=0$, $tT(t)=-\gamma t^2$ for $t\geq0$ and $tT(t)=t^2\geq-\gamma t^2$ for $t<0$, so $T$ is star $\gamma$-hypomonotone at $x^\ast=0$. As in
Remark~\ref{rem:weak-minty}, $T$ satisfies
\hyperref[Z1]{\textnormal{(Z1)}}--\hyperref[Z4]{\textnormal{(Z4)}}.
However, for every $c\geq1/\gamma$ we have $1-c\gamma\leq0$, and hence
\[
 \ran(\Id+cT)=\{(1-c\gamma)t: \, t\geq0\}\cup\{(1+c)t: \, t<0\} =\,]-\infty,0] \neq \R.
\]
Moreover, Corollary~\ref{cor:weak-star} concerns existence only: for $T=-\gamma\Id$ on $\R$, which satisfies its hypotheses with $x^\ast=0$,
and $c_k\equiv c\in\,]0,1/\gamma[$, the iterates are $x_k=(1-c\gamma)^{-k}x_0$, which diverge whenever $x_0\neq0$.
\end{remark}

\begin{remark}[Star cohypomonotone operators]
\label{rem:star-cohypo}
 For every prescribed $c>0$, star $\rho$-cohypomonotonicity alone does not guarantee surjectivity of $\Id+cT$ (Remark~\ref{rem:weak-minty}). If, in addition, $\norm{u}\leq L\norm{x-x^\ast}$ for all $(x,u)\in\gph(T)$, then $T$ is star $\rho L^2$-hypomonotone at $x^\ast$, and Corollary \ref{cor:weak-star} yields $\ran(\Id+cT) = H$ for every $c \in \, ]0, 1/(\rho L^2)[$ whenever one of
\hyperref[Z1]{\textnormal{(Z1)}}--\hyperref[Z4]{\textnormal{(Z4)}}
holds. This contrasts with the pairwise setting: if $T$ is maximally $\rho$-cohypomonotone, then $cT$ is maximally $(-\rho/c)$-comonotone in the sense of \cite{BMW}, and hence $\ran(\Id+cT)=H$ for every $c>\rho$ \cite{BMW}. Thus pairwise cohypomonotonicity yields surjectivity for large parameters, whereas its radial counterpart yields it, under linear growth, only for small parameters.
\end{remark}

\section{The Inexact Proximal Point Algorithm} \label{sec:05}

\subsection{The Algorithm and Fundamental Estimates}

We call the sequences generated by Algorithm~\ref{alg:ppa} \emph{inexact proximal sequences}; any admissible choice is allowed at each iteration, including when the resolvent is set-valued. All convergence results below hold for every such sequence, without any selection rule, since each proximal step is compared only with a zero of $T$, not with another proximal step. A uniform positive lower bound on the proximal parameters, when assumed, is denoted by $\underline{c}>0$.

\begin{algorithm}[ht]
\caption{Inexact proximal point algorithm for a star-monotone operator} \label{alg:ppa}
\textbf{Input:} $x_0\in H$, positive parameters $\{c_k\}_{k\geq0}$, and nonnegative tolerances $\{\varepsilon_k\}_{k\geq0}$.

\textbf{Iteration:} Given $x_k$, choose $x_{k+1}\in H$ and $\xi_{k+1}\in T(x_{k+1})$ such that
\begin{equation}\label{iterates}
 \norm{x_{k+1}+c_k\xi_{k+1}-x_k}\leq\varepsilon_k.
\end{equation}
Set $e_k:=x_{k+1}+c_k\xi_{k+1}-x_k$. Increase $k$ and repeat.
\end{algorithm}
By \eqref{eq:resolvent}, the iterates generated by Algorithm~\ref{alg:ppa} satisfy
\begin{equation}\label{eq:ppa}
 x_{k+1}\in J_{c_k}(x_k+e_k),\qquad \norm{e_k}\leq\varepsilon_k.
\end{equation}
Equivalently, with the selected residuals $\xi_{k+1}\in T(x_{k+1})$,
\begin{equation}\label{eq:implicit}
 x_k+e_k=x_{k+1}+c_k\xi_{k+1},\qquad
 \norm{e_k}\leq\varepsilon_k.
\end{equation}
Zero errors recover the exact proximal point method. With
summable tolerances, \eqref{iterates} is the selected-residual form of
criterion~$(A')$ in \cite[p.~880]{rock-1976}.

The criterion~\eqref{iterates} can be checked a posteriori: for a candidate pair $(x_{k+1},\xi_{k+1})$ with $\xi_{k+1}\in T(x_{k+1})$, the error $e_k$ is computable, so $\norm{e_k}\leq\varepsilon_k$ can serve as a stopping rule for an inner solver.

For the estimates on the sum of squared residuals and the general
star-monotone convergence results, we assume summable error tolerances:
\begin{equation}\label{eq:summable-errors}
 \sum_{k=0}^{\infty} \varepsilon_k <+\infty.
\end{equation}
The strong star-monotone case also admits the weaker error conditions specified in Theorem~\ref{thm:strong-convergence}.
If $T$ is star-monotone at some point and satisfies one of \hyperref[Z1]{\textnormal{(Z1)}}--\hyperref[Z4]{\textnormal{(Z4)}}, then, by Theorem \ref{thm:surj-unified}, $J_{c_k}(x_k+e_k)\neq \emptyset$ for every perturbed point and every $c_k>0$, so Algorithm~\ref{alg:ppa} is well defined for every initial point and every positive parameter sequence.

If $T$ is maximally monotone, then $J_c$ is single-valued and firmly nonexpansive on $H$ \cite[Theorems 1.2--1.3]{C1}. For star-monotone operators, this pairwise property is replaced by a radial one with respect to $x^\ast$, which holds for every point of the possibly set-valued $J_c(p)$.

\begin{lemma}[Radial resolvent estimate] \label{prop:pointwise}
 Let $T$ be (strongly) star-monotone at $x^\ast$ with modulus $\gamma \geq 0$, and let $c>0$ and $p \in \ran(\Id+cT)$. Then every $z \in J_c(p)$ satisfies
\begin{equation} \label{eq:star-firm}
 \ip{p-x^\ast}{z-x^\ast} \geq (1+c\gamma)\norm{z-x^\ast}^2.
\end{equation}
In particular,
\begin{equation} \label{eq:radial-resolvent}
 \norm{z-x^\ast}\leq\frac{1}{1+c\gamma}\norm{p-x^\ast}.
\end{equation}
\end{lemma}

\begin{proof}
Since $u:=(p-z)/c\in T(z)$, star monotonicity gives $\ip{p-z}{z-x^\ast}\geq c\gamma\norm{z-x^\ast}^2$. Writing $p-z=(p-x^\ast)-(z-x^\ast)$ gives \eqref{eq:star-firm}, and the Cauchy--Schwarz inequality gives \eqref{eq:radial-resolvent}.
\end{proof}

For $\gamma=0$, \eqref{eq:star-firm} is equivalent to $\norm{z-x^\ast}^2+\norm{p-z}^2 \leq \norm{p-x^\ast}^2$, that is, every selection of $J_c$ is firmly quasi-nonexpansive with respect to $x^\ast$.

Strong star monotonicity alone does not recover the classical properties of monotone resolvents.

\begin{example}[Failure of classical resolvent properties] \label{rem:resolvent-failure}
Fix $\gamma \in \, ]0, 2[$ and define $T: \R \to \R$ by
\[
T(x):=
\begin{cases}
(\gamma+2)x-x^2, & x\leq2,\\[1mm]
\gamma x, & x>2.
\end{cases}
\]
Then $T$ is continuous and strongly star-monotone at $0$ with modulus $\gamma$, since $xT(x)-\gamma x^2=x^2(2-x)^+$. If $c > 1/(2-\gamma)$, $\Id+cT$ is not injective on $]-\infty,2]$, so $J_c$ is not single-valued; for instance, $J_2(6)=\{\frac32, 2\}$ when $\gamma=1$. For $\gamma=c=1$, $J_1$ is single-valued, but $p_k := 4 - k^{-2}$ gives $J_1(p_k)= 2-\frac1k$ and $J_1(4)=2$, so $|J_1(p_k)-J_1(4)|/|p_k - 4| = k \to +\infty$, and $J_1$ is not Lipschitz continuous.
\end{example}

If $T$ is $\rho$-hypomonotone and $c\rho<1$, then $J_c$ is single-valued and $(1-c\rho)^{-1}$-Lipschitz continuous on its domain \cite{IPS2003,Pennanen2002}. The operator in Example~\ref{rem:resolvent-failure} is $(2-\gamma)$-hypomonotone, and $\gamma=c=1$ is exactly the threshold $c\rho=1$. The analysis below compares each proximal step only with a zero of $T$, through Lemma~\ref{prop:pointwise}, and does not rely on estimates between two values of the resolvent.

\begin{proposition} \label{prop:fejer}
 Let $\{x_k\}_{k}$ be any sequence generated by Algorithm \ref{alg:ppa}, with errors $e_k$ and residuals $\xi_{k+1}$ satisfying \eqref{eq:implicit}. Then the following assertions hold:
\begin{enumerate}[label=\textnormal{(\alph*)},leftmargin=2.2em]
\item For every $z\in H$ and every $k\geq0$,
\begin{equation}\label{eq:fejer-id}
\begin{aligned}
 \norm{x_k+e_k-z}^2={}&\norm{x_{k+1}-z}^2
 +\norm{x_k+e_k-x_{k+1}}^2 +2c_k\ip{\xi_{k+1}}{x_{k+1}-z}.
\end{aligned}
\end{equation}
\item If $T$ is star-monotone at $z$ with modulus
$\gamma\geq0$, then
\begin{equation}\label{eq:fejer-estimate}
 (1+2\gamma c_k)\norm{x_{k+1}-z}^2
 +\norm{x_k+e_k-x_{k+1}}^2
 \leq\norm{x_k+e_k-z}^2.
\end{equation}
In particular,
\begin{equation}\label{eq:quasi-fejer-distance}
 \norm{x_{k+1}-z}\leq\norm{x_k-z}+\varepsilon_k.
\end{equation}
\end{enumerate}
\end{proposition}

\begin{proof}
Expand $x_k+e_k-z=(x_{k+1}-z)+c_k\xi_{k+1}$ to obtain
\eqref{eq:fejer-id}. Substitution of
$\ip{\xi_{k+1}}{x_{k+1}-z}\geq\gamma\norm{x_{k+1}-z}^2$
gives \eqref{eq:fejer-estimate}. Dropping its nonnegative terms and
using $\norm{x_k+e_k-z}\leq\norm{x_k-z}+\varepsilon_k$ proves
\eqref{eq:quasi-fejer-distance}.
\end{proof}

The preceding estimate allows summable increases in the distances to a zero. We use the corresponding squared-norm notion of quasi-Fej\'er monotonicity.

\begin{definition}[Quasi-Fej\'er monotonicity of order two; {\cite[Definition 1]{AIS1998}}]\label{def:quasi-fejer}
Let $C\subset H$ be nonempty. A sequence $\{y_k\}_{k} \subset H$ is \emph{quasi-Fej\'er monotone of order two with respect to $C$} if,
for every $z \in C$, there exist an integer $k_z \geq 0$ and a nonnegative sequence $\{\eta_k(z)\}_{k}$ such that
\[
 \sum_{k=0}^{\infty} \eta_k(z)<+\infty, \qquad  \norm{y_{k+1}-z}^2\leq\norm{y_k-z}^2+\eta_k(z) \quad(k\geq k_z).
\]
\end{definition}

\begin{proposition}[Quasi-Fej\'er behavior]
\label{prop:residual-budget}
 Let $C\subset\zer(T)$ be nonempty, and suppose that $T$ is star-monotone at every point of $C$. Let $\{x_k\}_{k}$ be any sequence generated by Algorithm \ref{alg:ppa}, satisfying
\eqref{eq:summable-errors}. For each $z\in C$, set
\[
 M_z := \norm{x_0-z} + \sum_{k=0}^{\infty} \varepsilon_k.
\]
Then $\{x_k\}_{k}$ is quasi-Fej\'er monotone of order two with respect to $C$, is bounded, and $\{\norm{x_k-z}\}_{k}$ converges for every $z\in C$.
More precisely, $\norm{x_k-z} \leq M_z$ and
\begin{equation}\label{eq:quasi-fejer-squared}
 \norm{x_{k+1}-z}^2 \leq \norm{x_k-z}^2 + 2 M_z \varepsilon_k + \varepsilon_k^2.
\end{equation}
\end{proposition}

\begin{proof}
Fix $z\in C$ and put $d_k:=\norm{x_k-z}$.
Iterating \eqref{eq:quasi-fejer-distance} gives
$d_k\leq d_0+\sum_{j=0}^{k-1}\varepsilon_j\leq M_z$.
Moreover,
\[
 d_{k+1}+\sum_{j=k+1}^{\infty} \varepsilon_j
 \leq d_k+\sum_{j=k}^{\infty} \varepsilon_j.
\]
This nonnegative sequence converges. Since the error tails tend to zero, $d_k$ converges. Squaring $d_{k+1}\leq d_k+\varepsilon_k$ gives
\eqref{eq:quasi-fejer-squared}, whose error terms are summable. This proves the asserted quasi-Fej\'er property.
\end{proof}

The proximal identity also controls the residuals
$\xi_{k+1}\in T(x_{k+1})$. The next estimate will identify cluster
points as zeros.

\begin{corollary}[Estimates for the residuals]
\label{cor:residual-estimates}
Under the assumptions of Proposition~\ref{prop:residual-budget}, let
$\xi_{k+1}$ be as in \eqref{eq:implicit}. For every $z\in C$ and $N\geq1$,
\begin{equation}\label{eq:finite-residual-budget}
 \norm{x_N-z}^2+
 \sum_{k=0}^{N-1}c_k^2\norm{\xi_{k+1}}^2
 \leq
 \left(\norm{x_0-z}+\sum_{k=0}^{N-1}\varepsilon_k\right)^2.
\end{equation}
Consequently, $\sum_{k=0}^{\infty} c_k^2 \norm{\xi_{k+1}}^2 \leq M_z^2$, $x_{k+1}-x_k\longrightarrow0$, and
\begin{equation}\label{eq:best-residual}
 \min_{0\leq k<N} \norm{\xi_{k+1}} \leq \frac{M_z}{\left(\sum_{k=0}^{N-1} c_k^2 \right)^{1/2}}.
\end{equation}
In particular, $c_k\geq\underline{c}>0$ gives
$\norm{\xi_{k+1}}\to0$ and the bound
$M_z/(\underline{c}\sqrt{N})$ in \eqref{eq:best-residual}.
More generally, $\sum_k c_k^2= +\infty$ implies
$\liminf_k\norm{\xi_{k+1}}=0$.
\end{corollary}

\begin{proof}
Fix $z\in C$ and set $d_k:=\norm{x_k-z}$ and
$P_k:=\sum_{j=0}^{k-1}\varepsilon_j$, with $P_0=0$.
Proposition~\ref{prop:fejer}(b), with $\gamma=0$, gives
\begin{equation*}
 d_{k+1}^2+c_k^2\norm{\xi_{k+1}}^2
 \leq(d_k+\varepsilon_k)^2.
\end{equation*}
Since $d_k\leq d_0+P_k$, summing yields
\begin{align*}
 d_N^2+\sum_{k=0}^{N-1}c_k^2\norm{\xi_{k+1}}^2
 &\leq d_0^2+\sum_{k=0}^{N-1}
       \bigl(2(d_0+P_k)\varepsilon_k+\varepsilon_k^2\bigr)\\
 &=d_0^2+\sum_{k=0}^{N-1}
       \bigl((d_0+P_{k+1})^2-(d_0+P_k)^2\bigr)\\
 &=(d_0+P_N)^2.
\end{align*}
This proves \eqref{eq:finite-residual-budget} and the series bound. Thus $c_k\xi_{k+1}=x_k+e_k-x_{k+1} \to 0$, since $e_k \to 0$, also $x_{k+1}-x_k \to 0$. Furthermore,
\[
 \left(\min_{0 \leq k<N} \norm{\xi_{k+1}}^2 \right) \sum_{k=0}^{N-1}c_k^2 \leq\sum_{k=0}^{N-1} c_k^2 \norm{\xi_{k+1}}^2\leq M_z^2.
\]
The series bound gives $\xi_{k+1}\to0$ when $c_k \geq \underline{c}>0$. If instead $\sum_k c_k^2 = + \infty$ and $\liminf_k \norm{\xi_{k+1}} > 0$, the same series would diverge, a contradiction.
\end{proof}

\subsection{Convergence under Strong Star Monotonicity}

Strong star monotonicity yields a scalar contraction estimate for every choice of $x_{k+1} \in J_{c_k}(x_k+e_k)$. The estimate also permits nonsummable errors under a decay condition scaled by the proximal parameters.

\begin{theorem}[Strong convergence and perturbation bounds] \label{thm:strong-convergence}
 Let $\{x_k\}_{k}$ be any sequence generated by Algorithm \ref{alg:ppa} and suppose that $T$ is strongly star-monotone at $x^\ast$ with modulus $\gamma>0$. Set $q_k := (1 + c_k \gamma)^{-1}$. Then the following assertions hold.
\begin{enumerate}[label=\textnormal{(\alph*)},leftmargin=2.2em]
\item For every $k\geq0$,
\begin{equation}\label{eq:one-step}
 \norm{x_{k+1}-x^\ast} \leq q_k \bigl(\norm{x_k-x^\ast} + \varepsilon_k \bigr).
\end{equation}
Consequently, for every $k\geq1$,
\begin{equation}\label{eq:product}
 \norm{x_k-x^\ast}\leq{} \left( \prod_{i=0}^{k-1} q_i\right)\norm{x_0-x^\ast} +\sum_{j=0}^{k-1} \varepsilon_j\prod_{i=j}^{k-1}q_i.
\end{equation}
\item Suppose that $\sum_kc_k=+\infty$. Then $x_k\to x^\ast$ under either of the following error conditions:
\[
 \sum_k \varepsilon_k < +\infty, \qquad \text{or} \qquad \frac{\varepsilon_k}{c_k} \longrightarrow 0.
\]
For the parameter-scaled error condition, a further bound holds: if $L := \limsup_k \varepsilon_k/c_k < +\infty$, then
\begin{equation}\label{eq:strong-noise-floor}
 \limsup_{k \to +\infty} \norm{x_k-x^\ast} \leq \frac{L}{\gamma}.
\end{equation}
\item If $c_k \geq \underline{c}>0$, set $q := (1+\underline{c}\gamma)^{-1}\in\,]0,1[$. Then,
\begin{equation}\label{eq:inexact-rate}
 \norm{x_k-x^\ast} \leq q^k\norm{x_0-x^\ast} + \sum_{j=0}^{k-1} q^{k-j} \varepsilon_j, \qquad k \geq 0.
\end{equation}
The sum in \eqref{eq:inexact-rate} bounds the accumulated effect of the error tolerances. In particular, $\varepsilon_k\to0$ suffices for $x_k \to x^\ast$. For the exact method,
\begin{equation}\label{eq:linear-rate}
 \norm{x_k-x^\ast} \leq q^k\norm{x_0-x^\ast}.
\end{equation}
\end{enumerate}
\end{theorem}

\begin{proof}
 Lemma~\ref{prop:pointwise}, applied at $x_k + e_k$, gives \eqref{eq:one-step}, while induction gives \eqref{eq:product}. Put $a_k:= \norm{x_k-x^\ast}$. If $\sum_kc_k = +\infty$, then  $\sum_k \log(1+\gamma c_k)=+\infty$: either $\gamma c_k>1$ infinitely often, or eventually $\log(1+\gamma c_k) \geq (\log 2) \gamma c_k$. Thus $\prod_{i=N}^{k-1} q_i \to 0$ for every fixed $N$. If the errors are summable, split the sum in \eqref{eq:product} at $j=N$. Its finitely many terms with $j<N$ tend to zero, and the remaining terms are bounded by $\sum_{j=N}^{ \infty} \varepsilon_j$. Letting $k \to +\infty$ and then $N \to +\infty$ proves convergence.

For the second error condition, rewrite \eqref{eq:one-step} as
\[
 a_{k+1} \leq q_ka_k+(1-q_k)\frac{\varepsilon_k}{\gamma c_k}.
\]
If $L<+\infty$, then for every $\delta>0$ there is $N$ such that
$\varepsilon_k/(\gamma c_k)\leq L/\gamma+\delta$ for $k\geq N$. Iteration gives, for $k>N$,
\[
 a_k\leq \left(\prod_{i=N}^{k-1}q_i\right)a_N
 +\left(1-\prod_{i=N}^{k-1}q_i\right)(L/\gamma + \delta).
\]
Pass to the upper limit and then let $\delta \downarrow 0$ to obtain \eqref{eq:strong-noise-floor}. Taking $L=0$ proves the second convergence assertion. Finally, $q_k\leq q$ under the lower bound on $c_k$, so \eqref{eq:product} gives \eqref{eq:inexact-rate}.
The remaining conclusions follow from part \textnormal{(b)} and \eqref{eq:one-step} with zero errors.
\end{proof}

The estimate \eqref{eq:inexact-rate} gives explicit rates once a decay rule for the tolerances is specified. We call the convergence to $x^\ast$ \emph{R-linear} if there exist $C \geq 0$ and $\tau \in \, ]0,1[$ such that $\norm{x_k -x^\ast} \leq C \tau^k$ for all sufficiently large $k$.
 
\begin{remark}[Rates and sharpness]\label{rem:strong-rates} For the rate estimates in \textnormal{(i)} and their sharpness
examples in \textnormal{(ii)}, assume $c_k\geq\underline{c}>0$
and set $q:=(1+\underline{c}\gamma)^{-1}$.
\begin{enumerate}[label=\textnormal{(\roman*)},leftmargin=2.2em]
\item If $\varepsilon_k\leq E\rho^k$ with $E\geq0$ and $\rho\in\,]0,1[$, summing the geometric series in \eqref{eq:inexact-rate} gives
\begin{align*}
 & \norm{x_k-x^\ast} \leq q^k\norm{x_0- x^\ast} + Eq\,\frac{q^k-\rho^k}{q-\rho} \quad (q \neq \rho), \\
 & \norm{x_k-x^\ast}\leq q^k\norm{x_0-x^\ast} +Ekq^k \quad (q=\rho),
\end{align*}
so the convergence is R-linear, with any rate $\tau \in \, ]\max\{q, \rho\},1[$. If $\varepsilon_k \leq E(k+1)^{-p}$ with $p>0$, the inequality $(k+1)/(k-\ell+1) \leq \ell + 1$ for $1\leq\ell\leq k$ gives
\[
 \norm{x_k - x^\ast} \leq q^k\norm{x_0-x^\ast} + \frac{EC_{q,p}}{(k+1)^p}, \qquad C_{q,p} := \sum_{\ell=1}^{\infty}q^\ell(\ell+1)^p<+\infty,
\]
which allows nonsummable tolerances when $0<p\leq1$. Linear rates are also known for optimistic dual extrapolation and restarted high-order methods \cite[Proposition~1]{SZZJM2020}, \cite[Theorem~3.2]{LJ2025}, for single-valued operators in finite dimensions under Lipschitz continuity or higher-order smoothness; the estimates above hold for the implicit, possibly inexact, step, without Lipschitz continuity and for set-valued $T$.

 \item For $H=\R$, $T= \gamma\Id$, $x_0\geq0$ and $e_k =\varepsilon_k \geq0$, \eqref{eq:one-step} and \eqref{eq:product} are equalities. Hence $\sum_kc_k=+\infty$ cannot be weakened; for $x_0=0$ and $c_k\equiv\underline c$, the bounds in \textnormal{(i)} are attained for $e_k=E\rho^k$, so the factor $k$ for $q=\rho$ cannot be removed with base $q$, and $e_k = E(k+1)^{-p}$ gives $x_k\geq qEk^{-p}$ for $k \geq 1$, so the polynomial order is sharp.

 \item With $T= \gamma\Id$, $x_0=1$, $c_k=(k+1)^{-1}$ and $e_k =\gamma c_k$, we get $x_k=1$ for every $k$, although $\varepsilon_k \to0$ and $\sum_k c_k = +\infty$. Hence \eqref{eq:strong-noise-floor} is attained with $L= \gamma$, and vanishing errors do not suffice when $c_k \to0$.
\end{enumerate}
\end{remark}

\subsection{Weak Convergence under Star Monotonicity}

For $\gamma=0$, the radial estimate does not give a strict contraction. We combine quasi-Fej\'er monotonicity with identification of weak cluster points as zeros. For this identification, we use the following closedness property at the value zero.

\begin{definition} \label{def:demiclosed-zero}
 An operator $T: H \rightrightarrows H$ is sequentially demiclosed at zero if, for every $x_n\rightharpoonup x$ and every $u_n\to0$ with $u_n \in T(x_n)$, one has $0\in T(x)$.
\end{definition}

For a single-valued fixed-point mapping $U$, the analogous condition for $\Id-U$ is commonly called weak regularity; see \cite[Definition 3.1 and Remark~3.2]{CRZ2018}. Sequential weak-to-strong closedness of the entire graph implies
this condition, but the latter concerns only the limiting value zero. Sequential demiclosedness at zero also implies norm closedness of $\zer(T)$, by applying the definition to zero values. In finite dimensions, it is enough to require that $x_n\to x$ and $u_n\to0$, with $u_n\in T(x_n)$, imply $0\in T(x)$; in particular, the usual closed-graph condition suffices.

Condition $(M)$, already used in the surjectivity framework \textnormal{(Z4)}, also supplies the closedness needed to identify weak cluster points.

\begin{proposition}[Condition $(M)$ supplies demiclosedness] \label{prop:M-demiclosed}
 If $T$ satisfies condition $(M)$, then its graph is sequentially weak-to-strong closed. In particular, $T$ is sequentially demiclosed at zero.
\end{proposition}

\begin{proof}
 Let $x_n \rightharpoonup x$, $u_n \to u$, and $u_n \in T(x_n)$. The sequence $\{x_n\}_{n}$ is bounded, and hence $\ip{u_n}{x_n}-\ip{u}{x} = \ip{u_n-u}{x_n} + \ip{u}{x_n-x} \longrightarrow 0$. Condition $(M)$ gives $u \in T(x)$.
\end{proof}

The quasi-Fej\'er property provides boundedness and convergence of distances. The following principle identifies when the whole sequence converges weakly. The argument follows the quasi-Fej\'er approach of \cite{AIS1998} (see also \cite{Comb2001}).

\begin{lemma}[Quasi-Fej\'er convergence principle; {\cite[Proposition 1]{AIS1998}}] \label{lem:fejer-weak}
 Let $C \subset H$ be nonempty and $\{y_k\}_{k} \subset H$ be quasi-Fej\'er monotone of order two with respect to $C$. Then the following assertions hold:
\begin{enumerate} [label=\textnormal{(\alph*)},leftmargin=2.2em]
 \item $\{y_k\}_{k}$ is bounded;

 \item $\{\norm{y_k-z}^2\}_{k}$ converges for every $z \in C$;

 \item if every weak sequential cluster point of $\{y_k\}_{k}$ belongs to $C$, then $\{y_k \}_{k}$ converges weakly to a point of $C$.
\end{enumerate}
\end{lemma}

\begin{theorem}[Weak convergence]\label{thm:weak-convergence}
 Let $T: H \rightrightarrows H$ have a nonempty zero set, and $\{x_k\}_{k}$ be any sequence generated by Algorithm \ref{alg:ppa}, satisfying
\eqref{eq:summable-errors}. Assume that:
\begin{enumerate}[label=\textnormal{(\alph*)},leftmargin=2.2em]
 \item $T$ is sequentially demiclosed at zero;

 \item $T$ is star-monotone at every point of $\zer(T)$;

 \item $c_k \geq \underline{c}>0$ for every $k$.
\end{enumerate}
 Then $x_k \rightharpoonup \bar x$ for some $\bar x \in \zer(T)$. Moreover, $\xi_{k+1} \to 0$, and the sum and minimum estimates in Corollary \ref{cor:residual-estimates} hold at every zero.
\end{theorem}

\begin{proof}
 Proposition \ref{prop:residual-budget}, with $C=\zer(T)$, gives quasi-Fej\'er monotonicity of order two, boundedness, and convergence of $\norm{x_k-z}$ for every $z\in\zer(T)$. Corollary \ref{cor:residual-estimates} yields $x_{k+1} -x_k\to0$ and $\xi_{k+1}\to0$. If $x_{k_j} \rightharpoonup\bar x$, then $x_{k_j+1} \rightharpoonup\bar x$. Since $\xi_{k_j+1} \in T(x_{k_j+1})$ and $\xi_{k_j+1} \to 0$, sequential demiclosedness at zero implies $0 \in T(\bar x)$. Thus, every weak sequential cluster point belongs to $\zer(T)$,  and Lemma \ref{lem:fejer-weak}\textnormal{(c)} gives weak convergence to a zero.
\end{proof}

\begin{corollary}[Well-posed weak convergence]
\label{cor:weak-convergence-frameworks}
 Suppose that $T$ is star-mo\-no\-to\-ne at every point of its nonempty zero set, is sequentially demiclosed at zero, and satisfies one of the four surjectivity frameworks in Theorem~\ref{thm:surj-unified}. If $c_k\geq \underline{c} >0$ and \eqref{eq:summable-errors} holds, then the inexact proximal point algorithm is globally well defined and every sequence generated by Algorithm \ref{alg:ppa} converges weakly to a zero of $T$. In alternative \textnormal{(Z4)}, sequential demiclosedness at zero need not be assumed separately.
\end{corollary}

\begin{proof}
 Theorem~\ref{thm:surj-unified} guarantees a resolvent value at each perturbed point, and Theorem~\ref{thm:weak-convergence} gives weak convergence. In \textnormal{(Z4)}, condition $(M)$ implies sequential demiclosedness at zero by Proposition \ref{prop:M-demiclosed}.
\end{proof}

\begin{remark}[Strong convergence] \label{rem:compact-orbit-convergence}
 Let $T$ be star-monotone at every point of its nonempty zero set, and assume that $x_n\to x$, $u_n\to0$ and $u_n\in T(x_n)$ imply $0\in T(x)$; in particular, $\zer(T)$ is norm closed. Let $\{x_k\}_k$ be generated by Algorithm \ref{alg:ppa} with \eqref{eq:summable-errors}. Then:
\begin{enumerate}[label=\textnormal{(\roman*)},leftmargin=2.2em]
 \item $\{x_k\}_k$ converges strongly to a zero of $T$ if and only if  $\liminf_k {\rm dist} (x_k,\zer(T))=0$.
 
 \item If $\{x_k\}_k$ is relatively norm compact and $\sum_k c_k^2 = +\infty$, then $\{x_k\}_k$ converges strongly to a zero of $T$, even when $c_k \to 0$. In finite dimensions, relative norm compactness follows from the boundedness of the iterates.
\end{enumerate}
 Indeed, by \eqref{eq:quasi-fejer-distance}, ${\rm dist} (x_k,\zer(T)) + \sum_{j\geq k} \varepsilon_j$ is nonincreasing and $\norm{x_n-x_k} \leq 2\,{\rm dist}(x_k,\zer(T)) + \sum_{j \geq k} \varepsilon_j$ for $n\geq k$, which gives \textnormal{(i)} (see also \cite{Comb2001}). For \textnormal{(ii)}, Corollary \ref{cor:residual-estimates} gives $\xi_{k_j+1} \to 0$ along a subsequence, and a further subsequence satisfies $x_{k_j+1}\to\bar x\in \zer(T)$, so the condition in \textnormal{(i)} holds. For exact proximal sequences of maximally monotone operators, $\sum_k c_k^2 = + \infty$ already gives weak convergence in Hilbert spaces \cite{BL1978}, and it cannot be replaced by $\sum_k c_k = + \infty$ \cite[Remark~1, p.~331]{BL1978}.
\end{remark}

\subsection{Failure of Convergence under Star Monotonicity at a Single Zero}

The convergence result in Theorem~\ref{thm:weak-convergence} requires star monotonicity with respect to the whole zero set. In this subsection, we show that star monotonicity at a single zero cannot, in general, replace this assumption. We use a sequence generated by Algorithm~\ref{alg:ppa} with zero errors at every iteration.

We first give an operator that is star-monotone at exactly one of its zeros.

\begin{example}\label{ex:depends-zero}
 Let $R(x_1,x_2):=(-x_2,x_1)$ and define
\begin{equation}\label{eq:counter-operator}
 T(x):=(\norm{x}-1)Rx,\qquad x\in\R^2.
\end{equation}
 This operator is continuous, $\zer(T)=\{0\} \cup \Sph^1$ and $\ip{T(x)}{x}=0$, so it is star-monotone at $0$. For $z\in\Sph^1$ and $r>1$, the point $x=-rRz$ satisfies $\norm{x} = r$ and $Rx=rz$, hence $\ip{T(x)}{x-z}=-r(r-1)<0$. Thus $0$ is its only zero at which star monotonicity holds.
\end{example}

\begin{example}[Failure of convergence] \label{ex:failure}
 Let $T$ be given by \eqref{eq:counter-operator}, fix $c>0$ and choose $x_0$ with $r_0:=\norm{x_0}>1$. The function $\psi_c(r) := r\sqrt{1 + c^2(r-1)^2}$ is continuous and strictly increasing on $[1,+\infty[$, with $\psi_c(1)=1$. Given $x_k$ with $r_k:= \norm{x_k}>1$, we have $\psi_c(1) <r_k < \psi_c (r_k)$. Hence there is a unique $r_{k+1}\in \, ]1, r_k[$ such that
\begin{equation}\label{eq:radial}
 r_k=r_{k+1}\sqrt{1+c^2(r_{k+1}-1)^2}.
\end{equation}
Define
\begin{equation}\label{eq:defxk}
 x_{k+1}:=\bigl(\Id+c(r_{k+1}-1)R\bigr)^{-1}x_k.
\end{equation}
For $a \in\R$ and $y\in\R^2$, $\norm{(\Id + aR)y}^2 = (1+a^2) \norm{y}^2$ and $(\Id+aR)^{-1}=(\Id-aR)/(1+a^2)$, since $R^2=-\Id$. Together with \eqref{eq:radial}, this gives $\norm{ x_{k+1}} = r_{k+1}$ and $x_k= x_{k+1} + cT(x_{k+1})$. Thus \eqref{eq:defxk} generates an exact sequence of Algorithm~\ref{alg:ppa}. The decreasing radii satisfy $r_k \to r_\infty \geq 1$; passing to the limit in \eqref{eq:radial} gives $r_\infty=1$.

Put $s_k:=r_k-1>0$. Rationalizing \eqref{eq:radial} yields
\[
 0< s_k-s_{k+1} \leq A s_{k+1}^2, \qquad A := \frac{c^2(1+s_0)}2.
\]
Dividing by $s_ks_{k+1}$ gives $s_{k+1}^{-1} -s_k^{-1} \leq As_{k+1}/s_k \leq A$. Summing yields $s_k\geq(s_0^{-1} + Ak)^{-1}$, so $\sum_k s_k = + \infty$. If $x_k \to \bar x$, then $\norm{\bar x} =1$ and $\ip{Rx_{k+1}}{R\bar x}\to1$. Telescoping $x_k-x_{k+1} = cs_{k+1} Rx_{k+1}$ therefore gives, for $M>N$ large,
\[
 \ip{x_N-x_M}{R\bar x} \geq \frac c2 \sum_{k=N}^{M-1} s_{k+1} \longrightarrow+\infty
 \quad (M \to +\infty),
\]
contrary to boundedness. Hence, the sequence $\{x_k\}_{k}$ does not converge, even when ${\rm dist} (x_k, \zer(T)) = s_k \to 0$ and $\norm{T(x_k)} = r_k s_k \to 0$. Thus the distance criterion for strong convergence can fail when star monotonicity is assumed only at one zero.

For $0<c<2\sqrt2$, the function $\psi_c^2$ is strictly increasing on $[0,+\infty[$, since $(r-1)(2r-1)\geq-1/8$; as $\norm{y+cT(y)} = \psi_c (\norm{y})$, it follows that $J_c$ is single-valued on $\R^2$, so this nonconvergent sequence is uniquely determined.
\end{example}

\subsection{Scope of the Convergence Results}\label{subsec:scope}

The examples illustrate convergence outside generalized monotonicity conditions.

\begin{example} \label{ex:not-pseudo}
 Let $R(x_1,x_2) := (-x_2, x_1)$ and define $T(x) := x +\norm{x}^2Rx$ on $\R^2$. Since $\ip{T(x)}{x} = \norm{x}^2$, the operator is strongly star-monotone at its unique zero $0$, with modulus $1$. For $a=(1,0)$ and $b=(3,-1)$, we have $T(a)=(1,1)$, $T(b)=(13,29)$, with $\ip{T(a)}{b-a}=1>0$ and $\ip{T(b)}{b-a}=-3 < 0$. Then, $T$ is not quasimonotone and, in particular, not pseudomonotone. Continuity gives alternative~\textnormal{(Z3)} of Theorem~\ref{thm:surj-unified}. Under the error and parameter conditions of Theorem \ref{thm:strong-convergence}, every sequence generated by
Algorithm \ref{alg:ppa} converges to $0$.
\end{example}

\begin{example}[Radial sinusoidal rotation] \label{ex:not-hadamard}
Let $R(x_1,x_2):=(-x_2,x_1)$ and define $T(0):=0$ and
\[  
 T(x) := x+ \sqrt{\norm{x}} \sin(1/\norm{x}) Rx, \qquad x \in \R^2 \setminus \{0\}.
\]
For $r=\norm{x}>0$, $\norm{T(x)} =r \sqrt{1+ r\sin^2(1/r)} \to 0$ as $r\to0$, and $\ip{T(x)}{x} = \norm{x}^2$, since $Rx\perp x$. Thus $T$ is continuous on $\R^2$ and strongly star-monotone at its unique zero $0$, with modulus $1$. The following additional properties hold.

\medskip

\noindent $(i)$: Hypomonotonicity with modulus $\rho$ would imply $S(x) := (J_T(x)+ J_T(x)^T)/2\succeq -\rho I$ at every differentiability point, where $J_T$ denotes the Jacobian of $T$. This follows by applying the defining inequality to $x+th,x$ and letting $t\to0$. However, for $r_k:=(2k\pi)^{-1}$, $x_k:=(r_k,0)$ and $b_k := \sqrt{2k\pi}$,
\[
 J_T(x_k)=\begin{pmatrix}1&0\\-b_k& 1 \end{pmatrix},
 \qquad \lambda_{\min}(S(x_k))=1-b_k/2\to-\infty.
\]
Thus $T$ is not hypomonotone on any neighborhood of $0$.

\medskip

\noindent $(ii)$: Since $y_k:=T(x_k)=x_k$ and $\det J_T(x_k)=1$, the inverse function theorem gives local inverses $g_k$ with $ Dg_k (y_k) = \begin{pmatrix}1&0\\b_k&1\end{pmatrix}$. Their symmetric parts also have minimum eigenvalue $1-b_k/2\to-\infty$. The same differential test excludes local hypomonotonicity of $T^{-1}$ near $(0,0)$, hence local cohypomonotonicity of $T$. The inverse hypomonotonicity assumptions in \cite{IPS2003,CP,BMW} therefore fail.

\medskip

\noindent $(iii)$: No positive multiple $\lambda T$ is locally submonotone in the sense of \cite[Definition 2, equation~(10)]{LT-2025}, even when its values are restricted to a neighborhood of zero as in \cite[Assumption 2(d)]{LT-2025}. Indeed, the same differential test shows that, at differentiability points, this condition with violation $\tau\geq0$ requires $M_{\lambda,\tau}:=\lambda S+\tau (I+ \lambda J_T)^T(I+\lambda J_T) \succeq 0$, whereas $\det M_{\lambda,\tau}(x_k)= \bigl( \lambda +\tau(1+ \lambda)^2\bigr)^2-2k\pi \lambda^2 (\tau+\frac14) \to-\infty$, with $x_k\to0$ and $\lambda T(x_k)\to0$.

\medskip

Theorem \ref{thm:strong-convergence} gives $x_k\to0$ whenever $\sum_k c_k = +\infty$ and either $\sum_k\varepsilon_k < +\infty$ or
$\varepsilon_k/c_k \to0$. If $c_k \geq \underline c>0$, the rate estimates in Remark \ref{rem:strong-rates} apply.
\end{example}

\section{Applications to Nonsmooth Minimax Problems}\label{sec:06}

For a convex--concave function, Rockafellar's saddle operator is maximally monotone and its zeros are the saddle points \cite{C1,Rockafellar1970}. This section treats the nonconvex--nonconcave problems described in the introduction. Subsection~\ref{subsec:gen-convexity} relates the generalized convexity notions involved to star monotonicity of the Clarke subdifferential, and Subsection \ref{subsec:minimax} applies Sections \ref{sec:03} and~\ref{sec:05} to the Clarke saddle operator. Minimization is the special case without coupling and without maximization variable (Remark~\ref{rem:minimization}).


\subsection{Star Quasiconvexity, Quasar-convexity, and the Clarke Subdifferential} \label{subsec:gen-convexity}

We recall two generalized convexity notions anchored at a minimizer. The quasar-convexity condition below is the derivative-free form used in \cite{ABKLY,BLL}, which already appeared in  \cite[Lemma~10]{HSS}.

\begin{definition} \label{def:anchored}
 Let $h: \R^{n} \to \R$, $\overline{x} \in \argmin_{\R^{n}}h$ and $\gamma \geq 0$. We say that $h$ is:
\begin{enumerate}[label=\textnormal{(\alph*)},leftmargin=2.2em]
 \item $(\kappa, \gamma)$-strongly quasar-convex with respect to $\overline{x}$, where $\kappa \in \, ]0,1]$, if
\begin{equation}\label{eq:def-quasar}
 h(\lambda \overline{x}+(1-\lambda)x) \leq \kappa \lambda h(\overline{x})+(1-\kappa\lambda) h(x) - \lambda \left(1-\frac{\lambda}{2-\kappa}\right) \frac{\kappa \gamma}{2} \norm{x-\overline{x}}^2;
\end{equation}

 \item (strongly) star quasiconvex with respect to $\overline{x}$ with modulus $\gamma$ \cite{Khanh-Lara-1}, if
\begin{equation}\label{eq:def-star-qcx}
 h(\lambda \overline{x}+(1-\lambda)x) \leq h(x) - \lambda (1-\lambda) \frac{\gamma}{2} \norm{x-\overline{x}}^2,
\end{equation}
\end{enumerate}
for all $\lambda \in [0,1]$ and all $x \in \R^{n}$. When $\gamma=0$, we simply say that $h$ is $\kappa$-quasar-convex, respectively star quasiconvex, with respect to $\overline{x}$.
\end{definition}


 When $\kappa=1$, Definition \ref{def:anchored}\textnormal{(a)} reduces to (strong) star-convexity \cite{NesterovPolyak2006}. For $\gamma=0$, star quasiconvexity with respect to $\overline{x}$ means that $h$ is nondecreasing along every ray from $\overline{x}$ \cite{Khanh-Lara-1}. Hence it holds for every quasiconvex function with a minimizer and for additively separable sums of such functions,  which need not be quasiconvex \cite{Lara-10}; this includes Bergson--Samuelson social welfare functions with quasiconvex individual utilities \cite{Lara-10} and the Cauchy loss $\sum_i \log(1+x_i^2)$ of robust regression \cite{Khanh-Lara-2}.

The following lemma and proposition relate these notions to star monotonicity of the Clarke subdifferential $\partial^C h$ \cite{Clarke1990}.

\begin{lemma}[{\cite{Khanh-Lara-2}}] \label{lem:star-qcx-clarke}
 Let $h: \R^n\to\R$ be locally Lipschitz, $\overline{x} \in\argmin_{\R^n}h$ and $\gamma \geq0$. Then $h$ is star quasiconvex with respect to $\overline{x}$ with modulus $\gamma$ if and only if $\partial^C h$ is star-monotone at $\overline{x}$ with modulus $\gamma/2$, that is, $\ip{v}{x-\overline{x}}\geq\frac{\gamma}{2}\norm{x-\overline{x}}^2$ for all $x\in\R^n$ and $v\in\partial^C h(x)$.
In particular, if $\gamma>0$, then $\zer (\partial^C h) = \argmin_{\R^n} h = \{\overline{x}\}$.
\end{lemma}

\begin{proposition} \label{prop:clarke-qconv}
 Let $h:\R^n\to\R$ be locally Lipschitz and $(\kappa,\gamma)$-strongly quasar-convex with respect to $x^\ast\in\argmin_{\R^n}h$, where $\kappa\in\,]0,1]$ and $\gamma\geq0$, and set $\mu:=\kappa\gamma/(2-\kappa)$. Then the following assertions hold:
\begin{enumerate}[label=\textnormal{(\alph*)},leftmargin=2.2em]
 \item $h$ is strongly star quasiconvex with respect to $x^\ast$ with modulus $2\mu$; hence $\partial^C h$ is star-monotone at $x^\ast$ with modulus $\mu$.
 \item For every $x\in\R^n$ and every $v\in\partial^C h(x)$,
\begin{equation}\label{eq:clarke-gap}
 \ip{v}{x-x^\ast}\geq\kappa\bigl(h(x)-h(x^\ast)\bigr).
\end{equation}
 \item $\zer(\partial^C h)=\argmin_{\R^n}h$.
\end{enumerate}
\end{proposition}

\begin{proof}
\textnormal{(a)}: By \cite[Remark~4(d)]{ABKLY}, we have
\begin{equation}\label{eq:quasar-growth}
 h(x)-h(x^\ast)\geq \frac{\kappa\gamma}{2(2-\kappa)}\norm{x-x^\ast}^2.
\end{equation}
Write the right-hand side of \eqref{eq:def-quasar} as $h(x)-\kappa\lambda(h(x)-h(x^\ast))-\lambda(1-\frac{\lambda}{2-\kappa})\frac{\kappa\gamma}{2}\norm{x-x^\ast}^2$. Since $\lambda\geq\lambda(1-\lambda)$ and $1-\frac{\lambda}{2-\kappa}\geq1-\lambda$, \eqref{eq:quasar-growth} gives $h(\lambda x^\ast+(1-\lambda)x)\leq h(x)-\lambda(1-\lambda)\bigl(\frac{\kappa^2\gamma}{2(2-\kappa)}+\frac{\kappa\gamma}{2}\bigr)\norm{x-x^\ast}^2=h(x)-\lambda(1-\lambda)\mu\norm{x-x^\ast}^2$, which is \eqref{eq:def-star-qcx} with modulus $2\mu$. Lemma~\ref{lem:star-qcx-clarke} completes the proof.

\textnormal{(b)}: Since $h$ is $\kappa$-quasar-convex, \cite[Theorem~31]{ABKLY} gives $h(x^\ast)\geq h(x)+\frac1\kappa\ip{v}{x^\ast-x}$, which is \eqref{eq:clarke-gap}.

\textnormal{(c)}: This follows from \cite[Proposition~33]{ABKLY} and Clarke's necessary condition \cite[Chapter~2]{Clarke1990}.
\end{proof}

For $\gamma=0$, however, star quasiconvexity does not control the Clarke stationary points, in contrast with Proposition~\ref{prop:clarke-qconv}\textnormal{(c)}.

\begin{remark}[Star quasiconvexity with $\gamma=0$] \label{rem:sqc-critical}
Let $h(t):=1+(|t|-1)^3$ on $\R$. Since $h$ is nondecreasing in $|t|$, it is star quasiconvex with respect to $0$, where $\{0\}=\argmin_\R h$, and $\partial^C h$ is star-monotone at $0$ by Lemma~\ref{lem:star-qcx-clarke}. However, $h'(1)=0$ and $h(1)=1>0$, so $1\in\zer(\partial^C h)\setminus\argmin_\R h$, and $x_k\equiv1$ is an exact proximal sequence for every choice of $c_k>0$. By Proposition~\ref{prop:clarke-qconv}\textnormal{(c)}, $h$ is not $\kappa$-quasar-convex for any $\kappa\in\,]0,1]$.
\end{remark}

\subsection{Proximal Points for Nonsmooth Minimax Problems} \label{subsec:minimax}

Let $H_1$ and $H_2$ be real Hilbert spaces, let $A:H_1\to H_2$ be bounded and linear, and let $\phi_0:H_1\to\R$ and $\psi_0:H_2\to\R$ be locally Lipschitz. We consider the minimax problem with saddle function $\mathcal L$, to be minimized in $x$ and maximized in $y$, and its Clarke saddle operator on $H:=H_1\times H_2$:
\begin{align}
 \mathcal L(x,y)&:=\phi_0(x)-\psi_0(y)-\ip{Ax}{y}, \label{eq:saddle-function}\\
 T_{\mathcal L}(x,y)&:=\bigl(\partial^C \phi_0 (x)-A^\ast y\bigr) \times \bigl(\partial^C \psi_0 (y) +Ax\bigr). \label{eq:saddle-operator}
\end{align}
By the exact sum rule for continuously differentiable perturbations \cite[Chapter~2]{Clarke1990}, $\zer(T_{\mathcal L})$ consists of the pairs $(x,y)$ such that $x$ is Clarke stationary for $\mathcal L(\cdot,y)$ and $y$ is Clarke stationary for $-\mathcal L(x,\cdot)$. In particular, every saddle point $a=(x^\ast,y^\ast)$ of $\mathcal L$, that is, $\mathcal L(x^\ast,y)\leq\mathcal L(x^\ast,y^\ast)\leq\mathcal L(x,y^\ast)$ for all $(x,y)$, is a zero of $T_{\mathcal L}$. Given such a point, we use the partial functions and the gap function
\begin{equation}\label{eq:partial-functions}
 \phi(x) := \mathcal L(x,y^\ast)-\mathcal L(x^\ast, y^\ast), ~~ \psi(y):=\mathcal L(x^\ast, y^\ast)-\mathcal L(x^\ast,y), ~~ G(x,y) := \phi(x)+\psi(y).
\end{equation}
Thus $\phi$ and $\psi$ are nonnegative, $x^\ast\in\argmin\phi$, $y^\ast\in\argmin\psi$, and $G(x,y)=\mathcal L(x,y^\ast)-\mathcal L(x^\ast,y)\geq0$ is the gap of $\mathcal L$ relative to $a$. Since $\phi=\phi_0-\ip{A^\ast y^\ast}{\cdot}$ and $\psi=\psi_0+\ip{Ax^\ast}{\cdot}$ up to additive constants, the skew operator $K(w_1,w_2):=(-A^\ast w_2,Aw_1)$ gives
\begin{equation}\label{eq:saddle-decomposition}
 T_{\mathcal L}(z)=D(z)+K(z-a),\qquad D(x,y):=\partial^C\phi(x)\times\partial^C\psi(y),\qquad z\in H.
\end{equation}
Since $\ip{Kw}{w}=0$ for every $w$, \eqref{eq:saddle-decomposition} yields
\begin{equation}\label{eq:saddle-cancellation}
 \ip{u}{z-a}=\ip{d}{z-a},\qquad u=d+K(z-a)\in T_{\mathcal L}(z),\quad d\in D(z).
\end{equation}
Consequently, $T_{\mathcal L}$ is (strongly) star-monotone at $a$ with modulus $\gamma\geq0$ if and only if $D$ is. The operator $T_{\mathcal L}$ itself does not depend on $a$; the radial condition is imposed only on the partial functions $\mathcal L(\cdot,y^\ast)$ and $-\mathcal L(x^\ast,\cdot)$ at the saddle point. This condition is the weak Minty condition of \cite{Bohm,DDJ} with zero parameter. Unlike \cite{GTHG,PMC}, the results below require neither maximality nor Lipschitz continuity, and they hold for every $c>0$.

\begin{theorem}[Proximal points for nonsmooth minimax problems] \label{thm:saddle-ppa}
Let $a=(x^\ast,y^\ast)$ be a saddle point of $\mathcal L$, and suppose that $D$ in \eqref{eq:saddle-decomposition} is star-monotone at $a$ with modulus $\gamma\geq0$ and satisfies Condition~$(M)$. Then the following assertions hold.
\begin{enumerate}[label=\textnormal{(\alph*)},leftmargin=2.2em]
 \item $\ran(\Id+cT_{\mathcal L})=H$ for every $c>0$; equivalently, for every $(p,q)\in H$ there is $(x,y)\in H$ such that $x$ is a Clarke stationary point of $\mathcal L(\cdot,y)+\frac{1}{2c}\norm{\cdot-p}^2$ and $y$ is a Clarke stationary point of $-\mathcal L(x,\cdot)+\frac{1}{2c}\norm{\cdot-q}^2$. In particular, Algorithm~\ref{alg:ppa} with $T=T_{\mathcal L}$ is well defined for every initial point, every positive parameter sequence and every tolerance sequence.

 \item $T_{\mathcal L}$ is star-monotone at $a$ with modulus $\gamma$. If $\gamma>0$, then $\zer(T_{\mathcal L})=\{a\}$, and every sequence $\{z_k\}_k$ generated by Algorithm~\ref{alg:ppa} with $T=T_{\mathcal L}$ satisfies
\begin{equation}\label{eq:saddle-contraction}
 \norm{z_{k+1}-a}\leq\frac{\norm{z_k-a}+\varepsilon_k}{1+\gamma c_k}, \qquad \norm{z_{k+1}-a}\leq\frac{1}{\gamma}\norm{\xi_{k+1}},
\end{equation}
together with all conclusions of Theorem~\ref{thm:strong-convergence} and Remark~\ref{rem:strong-rates}; moreover, $G(z_k)\to0$ whenever $z_k\to a$.
\end{enumerate}
\end{theorem}

\begin{proof}
\textnormal{(a)}: By \cite[Propositions~2.1.2 and~2.1.5]{Clarke1990}, $D$ is locally bounded, has nonempty weakly compact convex values, and has a norm-to-weak closed graph; hence it is norm-to-weak upper semicontinuous. Since $K$ is bounded and linear, $T_{\mathcal L}$ in \eqref{eq:saddle-decomposition} inherits these properties. Condition $(M)$ passes from $D$ to $T_{\mathcal L}$: if $z_\alpha\rightharpoonup z$, $u_\alpha\in T_{\mathcal L}(z_\alpha)$, $u_\alpha\rightharpoonup u$ and $\limsup_\alpha\ip{u_\alpha}{z_\alpha}\leq\ip{u}{z}$, then $d_\alpha:=u_\alpha-K(z_\alpha-a)\rightharpoonup u-K(z-a)$ by weak continuity of $K$, and since $\ip{K(z_\alpha-a)}{z_\alpha}=\ip{K(z_\alpha-a)}{a}\to\ip{K(z-a)}{a}=\ip{K(z-a)}{z}$ by skewness, we get $\limsup_\alpha\ip{d_\alpha}{z_\alpha}\leq\ip{u-K(z-a)}{z}$; hence $u-K(z-a)\in D(z)$, that is, $u\in T_{\mathcal L}(z)$. By \eqref{eq:saddle-cancellation}, $T_{\mathcal L}$ is star-monotone at $a$, and alternative \hyperref[Z4]{\textnormal{(Z4)}} of Theorem~\ref{thm:surj-unified} gives $\ran(\Id+cT_{\mathcal L})=H$ for every $c>0$. The equivalence with partial Clarke stationarity follows from the exact sum rule for continuously differentiable perturbations \cite[Chapter~2]{Clarke1990}, since $\partial^C_x\mathcal L(x,y)=\partial^C\phi_0(x)-A^\ast y$ and $\partial^C_y(-\mathcal L)(x,y)=\partial^C\psi_0(y)+Ax$. Applying this at the perturbed points $z_k+e_k$ shows that Algorithm~\ref{alg:ppa} is well defined.

\textnormal{(b)}: Since $K0=0$ and $0\in D(a)$, we have $0\in T_{\mathcal L}(a)$, and \eqref{eq:saddle-cancellation} transfers the star monotonicity of $D$ at $a$ to $T_{\mathcal L}$. If $\gamma>0$ and $0\in T_{\mathcal L}(z)$, then $0\geq\gamma\norm{z-a}^2$, so $z=a$. Theorem~\ref{thm:strong-convergence} and Remark~\ref{rem:strong-rates} now apply with $T=T_{\mathcal L}$ and $x^\ast=a$, which gives the first estimate in \eqref{eq:saddle-contraction} and the convergence assertions. The second estimate follows from $\gamma\norm{z_{k+1}-a}^2\leq\ip{\xi_{k+1}}{z_{k+1}-a}\leq\norm{\xi_{k+1}}\norm{z_{k+1}-a}$. Finally, $G$ is continuous and $G(a)=0$.
\end{proof}

The second estimate in \eqref{eq:saddle-contraction}, a global metric subregularity property \cite{Khanh-Lara-2}, is an a posteriori bound, since $\xi_{k+1}=(z_k+e_k-z_{k+1})/c_k$ is computable.

In finite dimensions, the hypotheses of Theorem~\ref{thm:saddle-ppa} reduce to conditions on the two partial functions.

\begin{corollary}[Star quasiconvex--concave minimax problems] \label{cor:saddle-sqc}
Let $H_1=\R^m$, $H_2=\R^n$, and let $a=(x^\ast,y^\ast)$ be a saddle point of $\mathcal L$. Suppose that $\mathcal L(\cdot,y^\ast)$ is star quasiconvex with respect to $x^\ast$ with modulus $\gamma_1\geq0$, and that $-\mathcal L(x^\ast,\cdot)$ is star quasiconvex with respect to $y^\ast$ with modulus $\gamma_2\geq0$. Then the conclusions of Theorem~\ref{thm:saddle-ppa} hold with $\gamma:=\min\{\gamma_1,\gamma_2\}/2$.
\end{corollary}

\begin{proof}
 The functions $\phi$ and $\psi$ in \eqref{eq:partial-functions} differ from $\mathcal L(\cdot,y^\ast)$ and $-\mathcal L(x^\ast, \cdot)$ by constants. Lemma \ref{lem:star-qcx-clarke} shows that $\partial^C \phi$ and $\partial^C\psi$ are star-monotone at $x^\ast$ and $y^\ast$ with moduli $\gamma_1/2$ and $\gamma_2/2$. Adding the two inequalities shows that $D$ is star-monotone at $a$ with modulus $\gamma$ (see also \cite{Khanh-Lara-2} for separable operators). Since $\gph(D)$ is closed \cite[Proposition~2.1.5]{Clarke1990}, Remark~\ref{rem:M-closed-graph} gives Condition~$(M)$.
\end{proof}

For $\gamma_1=\gamma_2=0$, Corollary~\ref{cor:saddle-sqc} guarantees only the existence of the proximal steps; Remark~\ref{rem:sqc-critical} shows that convergence to a saddle point cannot be expected in general. Quasar-convexity restores convergence and also provides a gap estimate.

\begin{corollary}[Quasar-convex--concave minimax problems] \label{cor:saddle-quasar}
Let $H_1=\R^m$, $H_2=\R^n$, and let $a=(x^\ast,y^\ast)$ be a saddle point of $\mathcal L$. Suppose that $\mathcal L(\cdot,y^\ast)$ is $\kappa_1$-quasar-convex with respect to $x^\ast$ and $-\mathcal L(x^\ast,\cdot)$ is $\kappa_2$-quasar-convex with respect to $y^\ast$, and set $\kappa:=\min\{\kappa_1,\kappa_2\}$. Then assertion \textnormal{(a)} of Theorem~\ref{thm:saddle-ppa} holds, and every sequence $\{z_k\}_k$ generated by Algorithm~\ref{alg:ppa} with $T=T_{\mathcal L}$ and $\sum_k\varepsilon_k<+\infty$ satisfies, with $M:=\norm{z_0-a}+\sum_k\varepsilon_k$,
\begin{equation}\label{eq:saddle-gap}
 \min_{0\leq k<N}G(z_{k+1})\leq\frac{M^2}{2\kappa\sum_{k=0}^{N-1}c_k},\qquad N\geq1.
\end{equation}
If, moreover, $\argmin_{\R^m}\mathcal L(\cdot,y^\ast)=\{x^\ast\}$, $\argmax_{\R^n}\mathcal L(x^\ast,\cdot)=\{y^\ast\}$ and $\sum_kc_k=+\infty$, then $z_k\to a$.
\end{corollary}

\begin{proof}
By Proposition~\ref{prop:clarke-qconv}\textnormal{(a)}, $D$ is star-monotone at $a$, and Condition~$(M)$ holds as in the proof of Corollary~\ref{cor:saddle-sqc}; this gives assertion \textnormal{(a)} of Theorem~\ref{thm:saddle-ppa}.

Write $z_{k+1}=(x_{k+1}, y_{k+1})$ and, by \eqref{eq:saddle-decomposition}, $\xi_{k+1}=d_{k+1}+K(z_{k+1}-a)$ with $d_{k+1}=(u,v)\in D(z_{k+1})$. Since $\phi, \psi\geq0$, the gap estimate \eqref{eq:clarke-gap}, applied to $\phi$ and $\psi$, and \eqref{eq:saddle-cancellation} give
\begin{align*}
 \ip{\xi_{k+1}}{z_{k+1}-a} & = \ip{u}{x_{k+1}-x^\ast} + \ip{v}{y_{k+1}-y^\ast} \\ 
 & \geq\kappa_1 \phi(x_{k+1}) + \kappa_2 \psi (y_{k+1}) \\
 & \geq \kappa G(z_{k+1}).
\end{align*}
Inserting this inequality in \eqref{eq:fejer-id} with $z=a$ and summing as in the proof of Corollary \ref{cor:residual-estimates} yields
\begin{equation}\label{eq:saddle-gap-budget}
 \norm{z_N-a}^2+\sum_{k=0}^{N-1} c_k^2 \norm{\xi_{k+1}}^2 + 2\kappa \sum_{k=0}^{N-1} c_k G(z_{k+1}) \leq \left(\norm{z_0-a} + \sum_{k=0}^{N-1} \varepsilon_k\right)^2.
\end{equation}
All terms on the left are nonnegative. Taking the smallest gap proves \eqref{eq:saddle-gap}. For the last assertion, since $T_{\mathcal L}$ is star-monotone at $a\in\zer(T_{\mathcal L})$, Proposition \ref{prop:residual-budget} with $C=\{a\}$ shows that $\{z_k\}_k$ is bounded and that $\norm{z_k-a}$ converges. By \eqref{eq:saddle-gap-budget}, $\sum_k c_kG(z_{k+1})<+\infty$. Since $G\geq0$ and $\sum_k c_k=+\infty$, we have $\liminf_{k\to\infty}G(z_{k+1})=0$. Thus there is a subsequence with $G(z_{k_j+1})\to0$. By boundedness, a further subsequence converges to some $\bar z=(\bar x,\bar y)$. Continuity gives $G(\bar z)=0$, that is, $\phi(\bar x) =\psi(\bar y)=0$. Hence $\bar x\in\argmin\mathcal L(\cdot, y^\ast)=\{x^\ast\}$ and $\bar y\in\argmax\mathcal L(x^\ast, \cdot)=\{y^\ast\}$. Thus $\bar z=a$, and the limit of $\norm{z_k-a}$ is zero.
\end{proof}

The uniqueness assumption is needed because, by \eqref{eq:saddle-cancellation}, the coupling cancels only against $a$; for minimization, see \cite{BLL}.

\begin{example}[A nonsmooth nonconvex--nonconcave game] \label{ex:nonconvex-game}
For $\alpha\in\R$, let $\vartheta(t):=|t|+t^2-t\cos t +\sin t$ and $\mathcal L_\alpha(x,y) := \vartheta(x)-\vartheta(y)-\alpha xy$ on $\R^2$, so that $\phi_0=\psi_0=\vartheta$ and $A(t) = \alpha t$. Then $\partial^C\vartheta(0)=[-1,1]$ and $\partial^C\vartheta(t)=\{\operatorname{sgn}(t)+t(2+\sin t)\}$ for $t\neq0$, so $tu\geq t^2$ for all $u\in\partial^C\vartheta(t)$. Hence $(0,0)$ is a saddle point, $\phi=\psi=\vartheta$, and $D$ is strongly star-monotone at $(0,0)$ with modulus $1$. Since $\vartheta''(t)=2+\sin t+t\cos t$ equals $2-(2k+1)\pi$ at $t=(2k+1)\pi$, it is unbounded below; hence $\mathcal L_\alpha$ is nonconvex--nonconcave, $T_{\mathcal L_\alpha}$ is not hypomonotone for any modulus (the symmetric part of its Jacobian at $(t,t)$ is $\vartheta''(t)\Id$, and the differential test of Example~\ref{ex:not-hadamard}\,$(i)$ applies), and $\vartheta'$ is not globally Lipschitz.
Nevertheless, by Corollary~\ref{cor:saddle-sqc}, for every $\alpha\in\R$ and $c>0$ all proximal steps exist, and every inexact proximal sequence satisfies $\norm{z_{k+1}} \leq (\norm{z_k} +\varepsilon_k)/(1+c_k)$.
\end{example}

\begin{remark}[Nonsmooth minimization] \label{rem:minimization}
 Let $h:\R^n\to\R$ be locally Lipschitz and strongly star quasiconvex with respect to $\overline{x}$ with modulus $\gamma>0$. For $H_2=\{0\}$ and $\phi_0=h$, we get $T_{\mathcal L}=\partial^C h$, which is strongly star-monotone at $\overline{x}$ with modulus $\gamma/2$ by Lemma~\ref{lem:star-qcx-clarke}; hence Theorem~\ref{thm:strong-convergence} and Remark~\ref{rem:strong-rates} apply. Every Clarke stationary point of the coercive function $y\mapsto h(y)+\frac{1}{2c_k}\norm{y-(x_k+e_k)}^2$, in particular every global minimizer, satisfies $x_k+e_k\in y+c_k\partial^C h(y)$, so it is an admissible step of Algorithm~\ref{alg:ppa}, and a global solution of the subproblems is unnecessary. This covers $(\kappa,\gamma')$-strongly quasar-convex functions, with $\gamma=2\kappa\gamma'/(2-\kappa)$ by Proposition~\ref{prop:clarke-qconv}\textnormal{(a)}, for which the exact method is studied in \cite[Theorems~18 and~21]{BLL} (see also \cite[Section~4]{ABKLY}); linear convergence of first-order methods for strongly star quasiconvex functions is established in \cite{Khanh-Lara-1}.
\end{remark}

The following example, motivated by robust representation learning, is locally Lipschitz, nonsmooth, nonconvex and not even star-convex, yet strongly star quasiconvex, so the inexact proximal point algorithm converges R-linearly while each subproblem is solved only to Clarke stationarity.

\begin{example}[Robust feature alignment] \label{ex:robust-alignment}
 Let $\tau,\beta>0$ and 
 $$\Phi(z):=\min\{\norm{z}, \tau\} +\beta\norm{z}^2,$$ 
 on $\R^d$, a stabilized capped loss used to limit the influence of outliers. Then $\Phi$ is locally Lipschitz, nonsmooth at $\norm{z}\in\{0,\tau\}$, and nonconvex, since its radial slope drops from $1+2\beta\tau$ to $2\beta\tau$ at $\norm{z}=\tau$; if $\beta\tau<1/2$, it is not even star-convex, since $\Phi(\tau u)>\frac12 \Phi(2\tau u)$ for $\norm{u}=1$. For $z\neq0$, every $v\in \partial^C \Phi(z)$ has the form $v= s\,z/\norm{z}$ with $s\geq2\beta\norm{z}$, so $\ip{v}{z} \geq2\beta\norm{z}^2$; hence $\partial^C \Phi$ is star-monotone at $0$ with modulus $2\beta$, and $\Phi$ is strongly star quasiconvex with modulus $4\beta$ by Lemma \ref{lem:star-qcx-clarke}.

In feature distillation, given student and teacher features $z_{s,i}\in\R^{d_s}$ and $z_{t,i}\in\R^{d_t}$, $i=1,\ldots,N$, one minimizes $h(W):=\frac1N\sum_{i=1}^N \Phi(Wz_{s,i}-z_{t,i})$ over $W\in\R^{d_t\times d_s}$. Assume that $Z_s:=[z_{s,1}, \ldots, z_{s,N}]$ has full row rank and that $W^\ast z_{s,i}=z_{t,i}$ for all $i$. By the Clarke sum and chain rules \cite[Chapter~2]{Clarke1990}, every $V\in\partial^C h(W)$ satisfies $V=\frac1N \sum_iv_iz_{s,i}^\top$ with $v_i\in\partial^C \Phi((W-W^\ast) z_{s,i})$, and thus
\begin{align*}
 \ip{V}{W-W^\ast} = \frac1N\sum_{i=1}^N \ip{v_i}{(W-W^\ast)z_{s,i}} & \geq \frac{2\beta}{N} \norm{(W-W^\ast)Z_s}_F^2 \\
 & \geq \frac{2\beta\,\sigma_{\min}(Z_s)^2}{N} \norm{W-W^\ast}_F^2.
\end{align*}
Hence, $\partial^C h$ is strongly star-monotone at $W^\ast$ with modulus $2\beta\,\sigma_{\min}(Z_s)^2/N$, and, by Remark \ref{rem:minimization}, every inexact proximal sequence, with subproblems solved only to Clarke stationarity, converges R-linearly to $W^\ast$ for constant parameters and geometrically decreasing tolerances.
\end{example}

\section{Conclusions} \label{sec:07}

Star monotonicity yields surjectivity of $\Id+cT$ for every $c>0$ under four incomparable topological frameworks, so every proximal step exists for possibly set-valued resolvents. Every inexact proximal sequence then converges weakly to a zero under star monotonicity at every zero, and converges strongly with explicit rates under strong star monotonicity, while one zero does not suffice. For nonsmooth nonconvex--nonconcave minimax problems, these results apply whenever the partial functions are star quasiconvex at a saddle point, with R-linear convergence in the strongly star quasiconvex case.

Extensions to star cohypomonotone operators and to splitting methods are left for future work.

\section*{Acknowledgments}

We thank R.T. Rockafellar and R. Kostecki for the hard copy of reference \cite{Zaran}. This work has been partially supported by ANID---Chile through Fondecyt Regular 1241040 (Lara).

\end{document}